\def\IR{{\mathbb R}}

\def\IA{{\mathbb A}} 
\def\IK{{\mathbb K}}

\def\IC{\mathbb C} 
\def\ID{{\mathbb D}} 
\def\IQ{{\mathbb Q}}

\def\zbar{{\overline{z}}}

\documentclass[12pt]{article} 

\usepackage[psamsfonts]{amssymb} 
\usepackage{amsfonts,amsmath}  

\usepackage{epsfig,multicol}

\newtheorem{theorem}{Theorem}
\newtheorem{lemma}{Lemma}
\newtheorem{corollary}{Corollary}
\usepackage{apptools}
\AtAppendix{\counterwithin{lemma}{section}}

\title{On the properness of $p$-conformal energy on the Teichm\"uller space of a Riemann surface.}
\author{Hala Alaqad, Jianhua Gong, Gaven Martin \& Cong Yao\thanks{
\newline
Work of the first two authors partially supported by the New Zealand Marsden Fund and U.A.E.U UPAR  grant (12S127). Yao was supported by the Fundamental Research Funds for the Central Universities, the Young Scientist Program of the Ministry of Science and Technology of China (No. 2021YFA1002200), the National Natural Science Foundation of China (No. 12401096, No. 12101362), the Natural Science Foundation of Shandong Province (No. ZR2024QA035, No. ZR2022YQ01). Martin was partially supported by the Humboldt Foundation and the Max Planck Institute, Pl\"on.
 \newline
 {\bf Statements and Declarations:} the authors have no competing interests as defined by Springer, or other interests that might be perceived to influence the results and/or discussion reported in this paper.
 \newline
{\bf Author addresses.}
JG: Department of Mathematics, UAE University, UAE.
\newline
email:  j.gong@uaeu.ac.ae
\newline
 GM:Institute for Advanced Study, 
Massey University,
New Zealand
\newline
email: g.j.martin@massey.ac.nz, 
Orcid:0000-0002-9395-5855
\newline
CY: Research Center for Mathematics and Interdisciplinary Sciences, Shan-
dong University, 266237, Qingdao and Frontiers Science Center for Nonlinear
Expectations, Ministry of Education, P. R. China
email: c.yao@sdu.edu.cn   
}
} 
\begin{document}
\maketitle 

\begin{abstract}  We establish that the $p$-conformal energy, $p\geq 1$, defined by the $L^p$-norms of the distortion of Sobolev mappings,  is a proper functional on the Teichm\"uller space of Riemann surfaces of a fixed genus. This result is an application of a result herein identifying explicitly both the unique extremal mappings of finite distortion between hyperbolic annuli of given modulus, and their $p$-conformal energy.
\end{abstract}
\newpage
\section{Introduction}  Let $\Sigma_i$, $i=1,2$, be two homeomorphic Riemann surfaces of finite genus $g\geq 2$.  If $f:\Sigma_1\to\Sigma_2$ is a monotone Sobolev mapping of finite distortion (see \S\ref{fd} below), the $p$-conformal energy of $f$, $p\geq 1$,  is
\begin{equation}
E_p[f]=\left[ \frac{1}{4\pi(g-1)} \int_{\Sigma_1} \IK^p(z,f)\; d\sigma_1 \right]^{1/p},\quad \IK(z,f)=\frac{|\partial_z f|^2+|\partial_\zbar f|^2}{|\partial_z f|^2+|\partial_\zbar f|^2}.
\end{equation}
Here $\IK(z,f)$ is the distortion of the Sobolev homeomorphism $f$ and of course $4\pi(g-1)$ is the hyperbolic area of $\Sigma_1$. These energies were first considered by L.V. Ahlfors in his celebrated proof of Teichm\"uller's theorem, \cite{Ahlfors}.   

In this setting Teichm\"uller's theorem provides a unique extremal  $\IK_0$-quasiconformal homeomorphism $f_0:\Sigma_1\to\Sigma_2$ in any given homotopy class of homeomorphism of constant distortion.  This mapping is a barrier in that if $[f_0]$ denotes the homotopy class of $f_0$,  then
\begin{equation}
\inf_{f\in [f_0]} E_p[f] \leq \left[\frac{1}{|\Sigma_1|} \int_{\Sigma_1} \IK^p(z,f_0)\; d\sigma_1 \right]^{1/p} =  \IK_0 \end{equation}
That  the extremal Teichm\"uller mapping has constant distortion yields the last equality here.

\medskip

Teichm\"uller space ${\mathcal T}_g$ is the space of conformal or complex structures on a topological surface $\Sigma_0$ of genus $g$, and two are equivalent if there is a conformal map between them which
is homotopic to the identity.  A point in ${\mathcal T}_g$ is Riemann surface, together with a homotopy class of
homeomorphisms from this fixed surface $\Sigma_0$. Thus all homeomorphisms between two
points in ${\mathcal T}_g$ are homotopic by definition, and Teichm\"uller's problem is to 
minimize $\IK = {\rm ess \;sup}_\Sigma \{ \IK(z,f)\}$ over all quasiconformal maps between the two given points.
The minimum is easily seen to exist from standard properties of quasiconformal mappings (equicontinuity and lower semi-continuity of the distortion, \cite{AIM}) and  $\frac{1}{2}\log(\IK)$ is the {\em Teichm\"uller distance} $d_{Teich}(\Sigma_1,\Sigma_2)$ - a much studied metric on the moduli space of a Riemann surface. It is not obvious, but nevertheless true (see \cite{MY1}),  that 
\[ \inf_{f\in [f_0]} E_p[f] \to \IK_0, \quad \mbox{ as $p\to \infty$}, \]
where $\IK_0$ is the distortion of the extremal quasiconformal mapping homotopic to $f_0$. However,  while Sobolev minimisers exist in a weak sense for $E_p$ in the homotopy class of $f_0$,  it is not known (but conjectured to be true) that they are homeomorphic, \cite{MY1}.

\medskip

Our interest in this article is in studying (with the obvious notation)
\begin{equation} \inf_{f\in [f_0:\Sigma_1\to\Sigma_2]} E_p[f]  \end{equation}
as either $\Sigma_1$ or $\Sigma_2$ varies its complex structure. Notice that in the case $p=\infty$ (the supremum norm) it is known that for any homotopy class
\[ \inf_{f\in [f_0:\Sigma_1\to\Sigma_2]} E_\infty[f] \to \infty \]
if  $\Sigma_1$ is fixed and $\Sigma_2$ tends to the boundary of moduli space, or vice versa.  This is the statement that the Teichm\"uller metric is complete. We ask if the energy $E_p[f]$ is also proper,  that is it grows to $\infty$ as one or other of $\Sigma_i$ tends to the boundary,  and seek provide a explicit growth growth estimates.

Notice that this is {\em not true} in the flat case. We recall  \cite[Corollary 4]{AIM2} which we adapt to these problems.
\begin{theorem} Let $\sigma>0$.  For any surjective mapping of finite distortion $f:\ID\setminus\{0\}\to \Omega$, $mod(\Omega)=\sigma$,   
\[ \int_{\ID\setminus\{0\}} \IK(z,f) \; dz\geq \pi \coth \frac{\sigma}{2\pi} \]
Up to rotation, there is a unique extremal homeomorphic mapping $f_0$ achieving the lower bound, and $f_0^{-1}$ is harmonic. 
\end{theorem}
However,  despite the fact that $\ID(0,\frac{1}{2})$ has finite area in the hyperbolic metric $d\sigma=\frac{1}{|z|^2\log^2(1/|z|)}$ of the punctures disk $\ID\setminus\{0\}$ we have
\[ \int_{\ID(0,\frac{1}{2})\setminus\{0\}} \IK(z,f) \; d\sigma =+\infty.\] 
Indeed, for the extremal mapping
\[ f_0(z)= \frac{z}{2|z|} (|z|+\sqrt{|z|^2+4}) :\ID\setminus \{0\}\to \{1<|w|< \frac{1}{2}(1+\sqrt{5}) \}\]
we cannot have finite conformal energy on any punctured sub-disk when integrated against $d\sigma$.  We calculate
\[ (f_0)_z= \frac{1}{2} (1+\sqrt{1+4/|z|^2})- \frac{1}{\sqrt{|z|^2+4}} \frac{1}{|z|}, \quad   
  (f_0)_\zbar = -\frac{1}{\sqrt{1+4/|z|^2}} \frac{1}{\bar z^2}   \]
\[ \IK(z,f)=\frac{|(f_0)_z|^2+|(f_0)_\zbar|^2}{|(f_0)_z|^2-|(f_0)_\zbar|^2} = \frac{|z|^2+2}{|z| \sqrt{|z|^2+4}}\approx \frac{2}{|z|} \in L^1(\ID)\]
and then note $\frac{2}{|z|^3\log^2(1/|z|)} \not\in L^1(\ID(0,\frac{1}{2}))$. In fact for all $1\leq p<\infty$ we have $\IK^p\in L^1(\ID)$.

\medskip

We address this problem by giving explicit bounds for the $p$-conformal energy of mappings of finite distortion between ``finite'' hyperbolic annuli.  We then  reduce to this case via the ``thick and thin'' decompositions of the domain surface.  We thus provide estimates on the growth of $E_p(\Sigma_1,\Sigma_2)$ and $E_p(\Sigma_2,\Sigma_1)$ for $\Sigma_1$ fixed and $\Sigma_2$ approaching the boundary through hyperbolic structures shrinking the length of a geodesic curve. We see a surprising large asymmetry in these numbers when $p=1$ which disappears as $p\to\infty$. 

\section{Statement of results.} Our main results are the following.

 \begin{theorem}\label{thm1}  Let $\Sigma_1, \Sigma_2$ be closed Riemann surfaces of genus $g$. Let $m = \frac{1}{2\pi} \max_A   mod(A)$ where $A$ is an annular region about a simple closed geodesic of $\Sigma_2$, 
 
\medskip 
\noindent{\bf (a) Short geodesic in domain.} Suppose that $\Sigma_1$ has a shortest geodesic of length 
\[  \ell \leq \frac{8}{(2m+1)e^{2(2m+1)} }\] 
If $f:\Sigma_1\to \Sigma_2$ is a homeomorphism of finite distortion,  then
\begin{equation}\label{mr1}
\int_{\Sigma_1} \IK(z,f) \; d\sigma(z) \geq    4\pi(g-1)-4+  \frac{\sqrt{2}}{2m+1}\Big(1+\log\frac{4}{\ell}\Big){\log  \frac{8}{\ell (2m+1)}}    \end{equation}

\medskip 

\noindent{\bf (b) Short geodesic in target.}
Suppose that $\Sigma_2$ has a shortest geodesic of length $\ell>0$. Then
\begin{equation}
\int_{\Sigma_1} \IK(z,f) d\sigma(z) \geq \frac{4\, a^2\, \pi}{2\ell}(\pi-\ell)
\end{equation}
\end{theorem}
And in the case of $L^p$ means of distortion we have the following.

 \begin{theorem}\label{thm3}  Let $\Sigma_1, \Sigma_2$ be closed Riemann surfaces of genus $g$. Let
$m = \frac{1}{2\pi} \max_A   mod(A)$ where $A$ is an annular region about a simple closed geodesic of $\Sigma_2$, and $\ell$ is the length of the shortest geodesic in $\Sigma_1$.\\
\noindent{\bf  Short geodesic in domain} Suppose that $\Sigma_1$ has a shortest geodesic of length 
 $ \ell \leq \frac{1}{10}.$ 
If $f:\Sigma_1\to \Sigma_2$ is a homeomorphism of finite distortion, 
\begin{equation}
\int_{\Sigma_1} \IK^{p}(z,f) d\sigma_1    \geq   \big(4\pi(g-1)-4\big) + \frac{\pi^p}{4^p}  \frac{ 2\ell^{1-p}}{( \pi m)^p} 
\end{equation}
\end{theorem}

\noindent {\bf Remarks.} 
\begin{enumerate} 
\item We expect that the estimates of Theorem \ref{thm1} (\ref{mr1}) are of the right order (up to $\log \log \frac{1}{\ell}$ terms.  Both the numbers $\ell$ and $m$ are uniformly bounded in a bounded region of moduli space.  In (a), $\Sigma_1$ will be ``near'' the boundary and in (b), $\Sigma_2$ will be ``near'' the boundary.  
\item In Theorem \ref{thm3} we see
\[ \|\IK(z,f)\|_{L^p(\Sigma_1)} \geq    \frac{  \ell^{(1-p)/p}}{ 2 m } \]
and thus $\|\IK(z,f)\|_{L^\infty(\Sigma_1)} \geq  \frac{ 1}{ 2 m \ell } $.
 \item We will see in Corollary \ref{cor2} a maximal collar about the shortest geodesic of $\Sigma_2$ has conformal modulus at least $\frac{4 \pi }{\ell}$ up to lower order terms. Thus the extremal quasiconformal mapping has $K\geq \frac{4 \pi }{\ell m}$.   
 \item The Teichm\"uller distance is at least the ratio of moduli and hence  $\approx \log \frac{4 \pi }{\ell m}$. Kerchoff's extensive investigation into the relationship between extremal modulus and the Teichm\"uller distance \cite{Kerchoff} shows this to be of the correct order. In turn, this shows the results of Theorem \ref{thm3} are of the correct order. 
 \item In the case that $\Sigma_2$ has the shortest geodesic,  the metric of $\Sigma_1$ is bounded above and below and so nothing new beyond the results of \cite{MM} in the Euclidean case occurs for $p>1$.

 \end{enumerate}

\section{Mappings of finite distortion.}\label{fd} Let $\Omega$ be a domain in $\IC$ and $f\in W^{1,1}_{loc}(\Omega)$,  the Sobolev space of mappings with locally integrable first derivatives. The mapping $f$ is said to have finite distortion on $\Omega$ if
\begin{enumerate}
\item $W^{1,1}_{loc}(\Omega)$.
\item $J(z,f)=|f_z|^2-|f_\zbar|^2 \in L^1_{loc}(\Omega)$; the Jacobian is locally integrable.
\item $\IK(z,f)=\frac{|f_z|^2+|f_\zbar|^2}{|f_z|^2-|f_\zbar|^2}$ is finite almost everywhere.
\end{enumerate} 
The mapping $f$ is said to be {\em quasiconformal} if it is homeomorphic and $\IK(z,f)\in L^\infty(\Omega)$. It is easy to see that this implies $f\in W^{1,2}_{loc}(\Omega)$ and rather more is true.  See \cite{AIM} for the modern theory of quasiconformal mappings and mappings of finite distortion.  Because of the conformal invariance of each of the requirements above, this 
 definition needs no modification so as to define mappings of finite distortion between Riemann surfaces.

\section{Annuli in Riemann surfaces.}

In this section we collect together a number of useful results we will need to make concrete estimates for the $p$-conformal energy of extremal mappings of finite distortion between Riemann surfaces.

\medskip

For $s>1$ let $A_s$ denote the {\em round} annulus
\begin{equation}
A_s=\big\{z\in\IC: 1/s< |z| <s \big\}.
\end{equation}
The hyperbolic metric on this annulus is given by the line element density
\begin{equation}\label{2} 
\lambda_s(z) = \frac{\pi}{2\log(s)}\; \frac{1}{|z| \cos\Big(\frac{\pi \log |z|}{2\log(s)}\Big)}
\end{equation} 
See Beardon and Minda, \cite[\S12.2]{BM}.  The length of the {\em core} geodesic $\{|z|=1\}$ is $\ell=2\pi \lambda_s(1)$, thus
\begin{equation}\label{seqn}
\ell = \frac{ \pi^2}{\log(s)}.
\end{equation}
As $\ell$ is a geometric invariant it is useful to write the density at (\ref{2}) as
\begin{equation} \label{4}
\lambda_{A_s}(z) = \frac{\ell}{2\pi }\; \frac{1}{|z| \cos\Big(\frac{\ell}{2\pi} \, \log |z| \Big)}.
\end{equation} 

For $1\leq r <s$ we can compute the hyperbolic distance 
\begin{eqnarray*}\label{rformula} \rho_{A_s}(1,r)& = &  \int_{1}^{r} \lambda_s(t)\, dt = \tanh ^{-1}\Big[\sin \big(\frac{\ell \log r}{2 \pi }\big)\Big] 
\end{eqnarray*}
 We put $\delta=\rho_{A_s}(1,r)$ and identify the inverse function as 
\begin{equation}\label{rformula2} r=\exp \left[\frac{2 \pi}{\ell} \sin^{-1}(\tanh (\delta )) \right] 
\end{equation}
Then with (\ref{circlen}) we find the hyperbolic length of any circle $\{|z|=r; 1/s<r<s \}$.  This can be formulated invariantly noting that doubly connected domains have a unique closed geodesic if the boundary is not a puncture ($s=\infty$) and are all conformal to a round annulus .  

\begin{lemma} Let $\Omega$ be a doubly connected domain in $\IC$ with boundary components not points.  Let $\ell$ be the length of the shortest geodesic $\alpha$.  Then $\{z\in \Omega : d_{hyp}(z,\alpha)=\delta\}$ is a Jordan curve of hyperbolic length 
\[  {\ell \, \cosh(\delta )} \]
\end{lemma}
\noindent{\bf Proof.} Using (\ref{4}) and (\ref{rformula2}) we compute the length of any circle $\{|z|=r\}$, $1/s< r <s$, in $A_s$ as 
\begin{equation}\label{circlen} 
2\pi r \times \frac{\ell}{2\pi }\; \frac{1}{r\cos\Big(\sin^{-1}(\tanh (\delta ))\Big)}  = {\ell \, \cosh(\delta )} \end{equation}
Next, $\Omega$ is isometric to $A_s$ with $s$ as at (\ref{seqn}).  This isometry is effected by a conformal mapping and so the image of the geodesic $\{|z|=r\}$ is a smooth Jordan curve. \hfill $\Box$
 
 \medskip
 
 The conformal modulus of the annulus $A(a,b)=\{a<|z|<b\}$ is defined to be 
 \[ mod(A) =  \log \frac{b}{a} \]
The uniformisation theorem implies that any annular Riemann surface (that is homeomorphic to a doubly connected domain in $\IC$ with non-degenerate boundary,  sometimes called a {\em ring}) is conformally equivalent to $A_s$ for the correct choice of $s$ and therefore has a well defined modulus, $2\log(s)$.  We want the following special case of this.

 \begin{theorem}\label{thm4} Let $A$ be a tubular neighbourhood of hyperbolic radius $\delta$ about a simple closed geodesic curve $\alpha$ of length $\ell$ in a Riemann surface $\Sigma$.  That is 
 \[ A= \{z\in \Sigma:\rho_{hyp}(z,\alpha)<\delta \} \]
 Suppose $\delta$ is small enough that $A$ is doubly connected. Then the conformal modulus of $A$ is
 \begin{equation}\label{11}
 \frac{4 \pi}{\ell} \sin^{-1}\big[ \tanh (\delta )   \big]
 \end{equation}
 and the area of $A$ is
 \begin{equation}\label{area}
 area(A) = 2 \ell \sinh \delta  
 \end{equation}
 \end{theorem}
 \noindent{\bf Proof.} By (\ref{seqn}) $A$ is conformal to $A_s$ with $\log(s)=\frac{\pi^2}{\ell}$.  We only need observe that $A_r\subset A_s$ is an isometric copy (and therefore has the same conformal modulus) when $r$ is chosen by (\ref{rformula}) and then   $mod(A_r)=2\log(r)$. The area calculation is well known and can most easily be seen by identifying the fundamental domain of a hyperbolic isometry,  giving quotient $A$, in the upper-half space.
 \hfill $\Box$
 
 \medskip
 Note that while the area of $A$ might be bounded,  the conformal modulus of $A$ can be arbitrarily large.
 We can usefully formulate the metric in a more invariant fashion using (\ref{seqn}) in (\ref{2}).
 
 \begin{lemma} Let $A$ be a doubly connected tubular neighbourhood about a simple closed geodesic curve $\alpha$ of length $\ell>0$ in a Riemann surface.  Then the hyperbolic metric density is given as 
 \begin{equation}\label{3}
\lambda_A(z) =\frac{\ell}{2\pi} \, \frac{\cosh (\delta )}{  \exp \left(\frac{2 \pi}{\ell}  \sin ^{-1}(\tanh (\delta )) \right)}, \quad d_{hyp}(z,\alpha) = \delta.
\end{equation}  
 \end{lemma}
 \noindent {\bf Proof.} Suppose first $A$ is a round annulus with central core geodesic the unit circle.  Then the metric density is given at (\ref{4}) above as
  \begin{equation}\label{4a}
\lambda_A(z) = \frac{\ell}{2\pi}\; \frac{1}{|z| \cos\Big(\frac{\ell \log |z|}{2\pi}\Big)}
\end{equation}
where we must have 
\[ -\frac{\pi}{2} < \frac{\ell \log |z|}{2\pi} < \frac{\pi}{2} \Rightarrow   e^{- \frac{\pi^2}{\ell}} <   |z|  < e^{\frac{\pi^2}{\ell}} \]
Substituting the formulas developed above,  in particular $r=|z|$ at (\ref{rformula2}),  give the invariant form described. \hfill $\Box$

 \section{The extremal mappings between annuli.}
 
 We now consider the annulus $A_r\subset A_s$ with area weight (from (\ref{4}))
 \begin{equation}
\omega(z)= \big(\frac{\ell}{2\pi}\big)^2 \frac{1}{|z|^2 \cos^2\Big(\frac{\ell \log |z|}{2 \pi}\Big)} 
 \end{equation} 
We fix $t>1$ and we wish to identify  the value
 \begin{equation}
 \inf_{f} \; \int_{A_r}  \IK^p(z,f) \, \omega(z) \; dz
 \end{equation}
where the infimum is taken over all homeomorphisms $f:A_r\to A_t$ of finite distortion (orientation preserving and mapping boundary to boundary).  

In \cite{MM} this problem is addressed and it is shown that the extremal mapping is radial, with an ordinary differential equation identifying the radial stretching of the mapping.   This is proved by identifying an equivalence between the Gr\"otzsch problem (mappings between rectangles) and the Nitsche problem (mappings of annuli) is discussed in \cite[\S2.3]{MM}. The easiest way forward is to use these results,  however we must conform to the normalisations there.  

Thus we set
\[ \IA_1= A(1,r^2), \quad \IA_2=A(1,t^2) \]
and define a weight (the incomplete hyperbolic metric on $\IA_1$) by
\begin{equation}\label{sig}
\sigma(w) =  \frac{\ell}{2\pi}\; \frac{1}{|w| \cos\Big(\frac{\ell \log(|w|/r)}{2 \pi}\Big)}, \quad w\in \IA_1.
\end{equation}
Since $rA_r=\IA_1$,  with $\sigma(w)|dw|$ defining a metric we find that $\IA_1$ is isometric to $A_r$ with the hyperbolic metric induced as a subset of $A_s$ with its complete metric. A simple check here is to see that the core geodesic is $\{|w|=r\}$ also with length $\ell$.

Now we must identify $\lambda$ from the equation \cite[(2.11)]{MM},
\begin{equation}\label{lam}
\sigma^2(w) = \frac{1}{4\pi^2}\; \lambda(z) e^{-4\pi \Re e(z)}, \quad e^{2\pi z}=w.
\end{equation}
We have $x = \Re e(z) = \frac{1}{2\pi}\log |w|$. Then from (\ref{sig}) and noting the slightly different exponential change of variables from \cite{MM}, we have 
\[   \; \frac{\ell^2}{|w|^2 \cos^2\Big(\frac{\ell \log(|w|/r)}{2 \pi}\Big)} =   \lambda(z) e^{-4\pi \Re e(z)} \]
Now we put $|w|=e^{2\pi x}$ and finally
\begin{equation}\label{weight}  \lambda(z)= \frac{\ell^2}{\cos^2\Big(\ell x -\frac{\ell  \log r}{2 \pi}\Big)}, \quad\quad x=\Re e(z). \end{equation}
Now $0<x<\frac{1}{\pi} \log(r)$ gives us
\begin{equation}\label{15}
-\frac{\pi}{2} < - \sin^{-1}[\tanh  \delta]   < \ell x-\sin^{-1}[\tanh  \delta  ] < \sin^{-1}[\tanh  \delta  ] <\frac{\pi}{2}
\end{equation}
as required.

\medskip

We check these formulae here as follows.  The identity will always be an extremal mapping with $\IK=1$.  Thus
we must have the hyperbolic area of $\IA_1$ equal to $\int_{\IQ_1}\lambda(x) dz$.  That is
\[
2\ell\sinh(\delta)=\ell^2\int_0^{\frac{1}{\pi}\log r}\frac{dx}{\cos^2\left(\ell x-\frac{\ell\log(r)}{2\pi}\right)}=2\ell \tan\left(\frac{\ell\log(r)}{2\pi}\right)
\]
as desired.

\medskip

After these changes and with
\[ \IQ_1=[0,\frac{1}{\pi}\log(r) ]\times[0,1], \quad \IQ_2=[0,\frac{1}{\pi}\log(t) ]\times[0,1],\]
 the problem becomes to identify
 \begin{equation}\label{24}
 \inf \; \int_{\IQ_1}  \IK^p(z,\tilde{f})\, \lambda(x) \; dz
 \end{equation}
where the infimum is taken over all homeomorphisms $\tilde{f}:\IQ_1\to \IQ_2$ of finite distortion mapping edges to edges and $\tilde{f}(0)=0$ . 

As shown in \cite{MM} the infimum is uniquely achieved by a homeomorphism of the form 
\begin{equation}\label{bv==} \tilde{f}(z) = u(x)+i y,  \quad u(0)=0,\; u\big(\frac{1}{\pi}\log r \big)=\frac{1}{\pi}\log(t) \end{equation}

At this point we see  (\ref{11}) of Theorem \ref{thm1} implies
\[ mod(\IA_1)=\frac{4 \pi}{\ell} \sin^{-1}\big[ \tanh (\delta )\big] \]
Thus we amend (\ref{weight}) to give to a more invariant formulation,
\begin{equation}\label{realweight}
\lambda(x) = \frac{\ell^2}{\cos^2\Big(\ell (x -\frac{1}{4\pi} mod(\IA_1))\Big)}, \quad  0\leq x \leq \frac{1}{2\pi} mod(\IA_1) \end{equation}
and with boundary values (\ref{bv==}),
\begin{equation}\label{bv=}
u(0)=0, \quad u\Big(\frac{1}{2\pi} mod(\IA_1)\Big) =  \frac{1}{2\pi} mod(\IA_2).
\end{equation}

\noindent{\bf Warning}.  Although this formulation is conformally invariant,  due  to the integration the minimum value of (\ref{24}) is not.  If $\varphi:\Omega \to \IA_2$ is a conformal mapping,  then 
\[ \IK(z,\varphi^{-1}\circ \tilde{f})=\IK(z, \tilde{f}) \]
and so the integral will not change if the target domain $\IA_2$ is replaced by any doubly connected domain of the same modulus.  However to obtain the form of the minimiser and the variational equation,  we must have that $\IA_1$ is a ``round'' annulus.

\medskip

\subsection{The form of minimisers} From \cite[Theorem 4.3]{MM}  $u$ is in fact  the solution to the equation
\begin{equation}\label{ve}
\lambda(x)\Big(1-\frac{1}{u_x^2}\Big)\Psi'(u_x+\frac{1}{u_x}) = \alpha, \quad \Psi(t)=t^p.
\end{equation}
 where $\alpha$ is a constant chosen so the boundary conditions (here surjectivity) are satisfied.  That is the boundary restrictions imposed by (\ref{bv=}). We note that as $\lambda(x)>0$ and $\Psi'(t)>0$, from (\ref{ve}) we obtain the following lemma.
 \begin{lemma} Exactly one of the following three situations must occur.
 \begin{enumerate}
 \item $u_x(x_0)=1$ for one value $x_0$, and hence all values $u_x(x)\equiv1$, and $\alpha=0$.
 \item $u_x(x_0)<1$ for one value $x_0$,  and hence all values $u_x(x)< 1$, and $\alpha\in (-\infty,0)$.
  \item $u_x(x_0)>1$ for one value $x_0$,  and hence all values $u_x(x)>1$, and $\alpha\in (0,\infty)$.
 \end{enumerate}
 \end{lemma} 
 Now (\ref{ve}) becomes with (\ref{realweight})
 \begin{equation}\label{veq}
\Big(1-\frac{1}{u_x^2}\Big)\Psi'(u_x+\frac{1}{u_x}) = \frac{\alpha}{\ell^2} \; \cos^2\Big(\ell (x -\frac{1}{4\pi} mod(\IA_1))\Big).
\end{equation}

\section{The case of minimisers of mean distortion.} 

In this special case where $p=1$ we can make more explicit calculations as $\Psi'(t)\equiv 1$ and so (\ref{veq}) reads as
\[ 
1-\frac{1}{u_x^2}= \frac{\alpha}{\ell^2} \cos^2\Big(\ell (x -\frac{1}{4\pi} mod(\IA_1))\Big),
\]
We are going to have to choose $\alpha$, given $\ell$, in what follows.  We have
\begin{equation}\label{uxeqn}
u_x = \frac{1}{\sqrt{1-\alpha \, \ell^{-2} \, \cos^2\Big(\ell (x -\frac{1}{4\pi} mod(\IA_1))\Big)}}
\end{equation}
Notice that this is not integrable on our range $x\in [0,\frac{1}{2\pi} mod(\IA_1))]$ if $\alpha\geq \ell^2$. Further,  we find that $u_x$ has a maximum where the weight has a minimum, a phenomenon  observed in \cite{MM} and that underpins the Nitsche conjecture,  see \cite{AIM2} for this in the flat case,  and related work by Kalaj \cite{Kalaj}.  We compute with 
\begin{equation}
\alpha_\ell = \frac{\alpha}{\ell^2}
\end{equation}
for notational simplicity.
\begin{eqnarray}
\int_0^{\frac{\ell}{2\pi} mod(\IA_1)} u_x(t)\;dt &=&   \frac{2}{\ell}  \int_{0}^{\frac{\ell}{4\pi} mod(\IA_1)} \frac{1}{\sqrt{1-\alpha_\ell \cos^2(s)}}  \;ds \nonumber \\
&=&\frac{2}{\ell}  \frac{F\left(\frac{\ell}{4\pi} mod(\IA_1) \left|\frac{\alpha_\ell }{1-\alpha_\ell}\right.\right)}{\sqrt{1-\alpha_\ell}} \label{boundaryvalues}
\end{eqnarray}
The result is in terms of the EllipticF function of the second kind.  We will later need the EllipticE function as well,  so we recall both of these.
\begin{eqnarray}\label{Fdef} F (\phi |m)=\int _0^{\phi }  \frac{d \theta}{\sqrt{1-m  \sin ^2\theta}},&&
\label{Edef} E (\phi |m)=\int _0^{\phi } {\sqrt{1-m  \sin ^2\theta}}\; d \theta
\end{eqnarray}
Near $0$ the integrand is $\frac{1}{\sqrt{1-\alpha_\ell+\alpha_\ell s^2/2}}$ and so the integral is strictly increasing on $\alpha_\ell\in (-\infty,1]$ and diverges to $+\infty$ as $\alpha_\ell \nearrow 1$, while the integral tends to $0$ as $\alpha_\ell \to -\infty$. Thus we can always solve uniquely for $\alpha_\ell$, and hence $\alpha$.

\begin{lemma}\label{lem3}  For each $\ell>0$ there is a unique $\alpha_\ell\in (-\infty,1)$ solving the boundary conditions given at (\ref{bv=}).
\end{lemma}

Performing the integration yields
\begin{equation}
u(x)=\frac{F\left(\ell (x-a)\left|\frac{\alpha_\ell }{\alpha_\ell -1}\right.\right)+F\left(a\ell\left|\frac{\alpha_\ell }{\alpha_\ell -1}\right.\right)}{\ell \sqrt{1-\alpha_\ell } }, \quad a=\frac{1}{4\pi} mod(\IA_1).
\end{equation}
If $m_\Omega=mod(\IA_2)$,  the conformal modulus of the target, we require $u (2a )=\frac{m_\Omega}{2\pi}$. That is 
\begin{equation}\label{alpha} \frac{2F\left(\frac{\ell}{4\pi} mod(\IA_1)\left|\frac{\alpha_\ell }{\alpha_\ell -1}\right.\right)}{\ell \sqrt{1-\alpha_\ell } } =\frac{m_\Omega}{2\pi}   \end{equation}

 \section{The integral of the distortion.}
We have $\IA_1$, a round hyperbolic annulus with core geodesic of length $\ell$ and hyperbolic radius $\delta$ in a Riemann surface $\Sigma_1$. We consider the extremal mapping $f:\IA_1\to \Sigma_2$ with image $f(\IA_1)=\IA_2$.  Due to conformal invariance of the integral the only invariant of relevance is $m_\Omega=mod(\IA_2)$. We have reduced the problem to a Gr\"otzch type problem for a mapping $\tilde{f}$ of explicit form which we have identified. 

Next we wish to calculate the integral of the distortion of $\tilde{f}$. We have
\begin{equation}
\IK(z,\tilde{f}) =\frac{1}{2}\left( u_x+\frac{1}{u_x} \right).
\end{equation}
With $T=\frac{1}{2\pi} mod(\IA_1)$ and $\lambda(s)$ as at (\ref{weight}) we want to compute
\begin{equation}\label{firstint}
\int_{\IQ_1} \IK(z,\tilde{f}) \lambda(z) dz = \frac{1}{2} \int_{0}^{T} \left( u_x(s)+\frac{1}{u_x(s)} \right) \lambda(s) ds.
\end{equation}
Set 
\begin{equation}\label{theta}
 \Theta = \frac{\ell T}{2} = \frac{\ell}{4\pi} mod(\IA_1)\end{equation}
We compute two integrals with (\ref{realweight}) and (\ref{uxeqn}) and the obvious change of variables.  
\begin{eqnarray}
\label{i1}I_1& = &  2\int_{0}^{T/2} u_x(s)\lambda(s) ds = 2 \ell^2
 \int_0^{T/2} \frac{\sec ^2\left(\ell x \right)}{\sqrt{1-\alpha_\ell  \cos ^2\left(\ell  x \right)}} \, dx\\
\nonumber&=&  {2 \ell} \left(\frac{F\left(\Theta \left|\frac{\alpha_\ell }{\alpha_\ell -1}\right.\right)}{\sqrt{1-\alpha_\ell }}-\sqrt{1-\alpha_\ell } E\left(\Theta \left|\frac{\alpha_\ell }{\alpha_\ell -1}\right.\right)+ \tan (\Theta ) \sqrt{1-\alpha_\ell  \cos ^2(\Theta )}\right)
\end{eqnarray}
Then 
\begin{eqnarray}
\label{i2} I_2& = &2\int_{0}^{T/2} \frac{1}{u_x(s)}\;\lambda(s) ds \\ \nonumber &=& {2}{\ell^2} \int_0^{T/2} \sec ^2(\ell x) \sqrt{1-\alpha_\ell  \cos ^2(\ell x )} \, dx \\
\nonumber & = & {2\ell} \Big( \sqrt{1-\alpha_\ell }  \Big[F\left(\Theta \left|\frac{\alpha_\ell }{\alpha_\ell -1}\right.\right)-E\left(\Theta \left|\frac{\alpha_\ell }{\alpha_\ell -1}\right.\right)\Big]+\tan (\Theta ) \sqrt{1-\alpha_\ell  \cos ^2(\Theta )}\Big)
\end{eqnarray}
Simplifying the sum $\frac{1}{2}(I_1+I_2)$ gives
\begin{eqnarray*}\frac{1}{2}(I_1+I_2) &=& {\ell} \left(\big(\sqrt{1-\alpha_\ell }+\frac{1}{\sqrt{1-\alpha_\ell }}\big) F\left[\Theta \left|\frac{\alpha_\ell }{\alpha_\ell -1}\right.\right]-2 \sqrt{1-\alpha_\ell } E\left[\Theta \left|\frac{\alpha_\ell }{\alpha_\ell -1}\right.\right]\right.\\ && \left.+2 \tan \Theta  \sqrt{1-\alpha_\ell  \cos ^2\Theta }\right)
\end{eqnarray*}
Next we recall that the choice of $\alpha_\ell$ at (\ref{alpha}) so that 
\begin{equation}\label{32} F\left(\Theta \left|\frac{\alpha_\ell }{\alpha_\ell -1}\right.\right) =\frac{m_\Omega}{4\pi}\ell \sqrt{1-\alpha_\ell }\end{equation}
\begin{eqnarray*}\frac{1}{2}(I_1+I_2) &=& {2\ell} \left((2-\alpha_\ell) \frac{m_\Omega\,\ell}{4\pi} -2 \sqrt{1-\alpha_\ell } E\left[\Theta \left|\frac{\alpha_\ell }{\alpha_\ell -1}\right.\right]\right.\\ && \left.+2 \tan \Theta  \sqrt{1-\alpha_\ell  \cos ^2\Theta }\right)
\end{eqnarray*}
\begin{center}
\scalebox{0.3}{\includegraphics{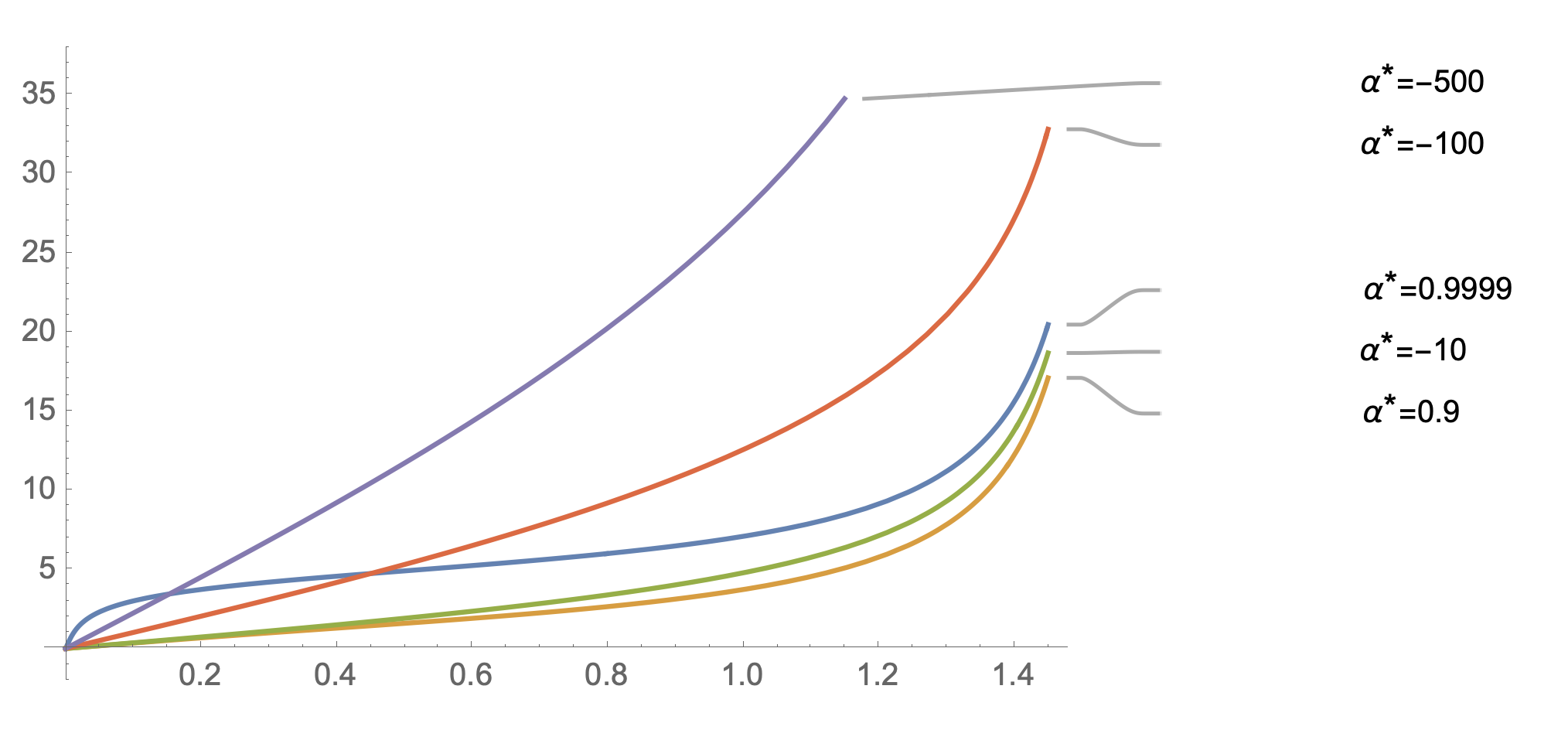}}
\end{center}
{\bf Figure.}{\em The graph of the function
\[\footnotesize{ \big(\sqrt{1-\alpha_\ell }+\frac{1}{\sqrt{1-\alpha_\ell }}\big) F\left[\Theta \left|\frac{\alpha_\ell }{\alpha_\ell -1}\right.\right]-2 \sqrt{1-\alpha_\ell } E\left[\Theta \left|\frac{\alpha_\ell }{\alpha_\ell -1}\right.\right] +2 \tan \Theta  \sqrt{1-\alpha_\ell  \cos ^2\Theta }}\]
for $\Theta\in [0,\frac{\pi}{2}]$ and a range of $\alpha_\ell$.}
 
 \bigskip

We have established the following theorem.

\begin{theorem}\label{thm2} Let $A$ be a tubular neighbourhood of hyperbolic radius $\delta$ about a simple closed geodesic curve $\gamma$ of length $\ell$ in a Riemann surface $\Sigma$.  That is 
 \[ A= \{z\in \Sigma:d_{hyp}(z,\gamma)<\delta \} \]
 Suppose that $A$ is doubly connected.  Let $\Omega$ be a doubly connected set of conformal modulus $m_\Omega$. If $f:A\to \Omega$ is a surjective proper mapping of finite distortion, then
 \begin{eqnarray}\nonumber
 \int_A \IK(z,f) & \geq& {\ell} \left(\big(\sqrt{1-\alpha_\ell }+\frac{1}{\sqrt{1-\alpha_\ell }}\big) F\big[\Theta \left|\frac{\alpha_\ell }{\alpha_\ell -1}\right.\big]  \right. \\
 && - 2 \sqrt{1-\alpha_\ell } E\big[\Theta \left|\frac{\alpha_\ell }{\alpha_\ell -1} \right.\big] \label{ineq} \\ \nonumber&& \left.+2 \tan \Theta  \sqrt{1-\alpha_\ell  \cos ^2\Theta }\right)  \end{eqnarray}
 where 
 \[ \Theta = \sin^{-1}(\tanh(\delta)) \in [0,\frac{\pi}{2}], \]
 and $\alpha_\ell$ is chosen so that
 \begin{equation}\label{alphaell}
 \frac{ F\left(\Theta\left|\frac{\alpha_\ell }{\alpha_\ell-1}\right.\right)}{\sqrt{1-\alpha_\ell } } = \frac{\ell \;m_\Omega}{4\pi}
 \end{equation}
 The inequality in (\ref{ineq}) is sharp,  it holds with equality for the unique extremal mapping.
\end{theorem}
\begin{center}
\scalebox{0.25}{\includegraphics{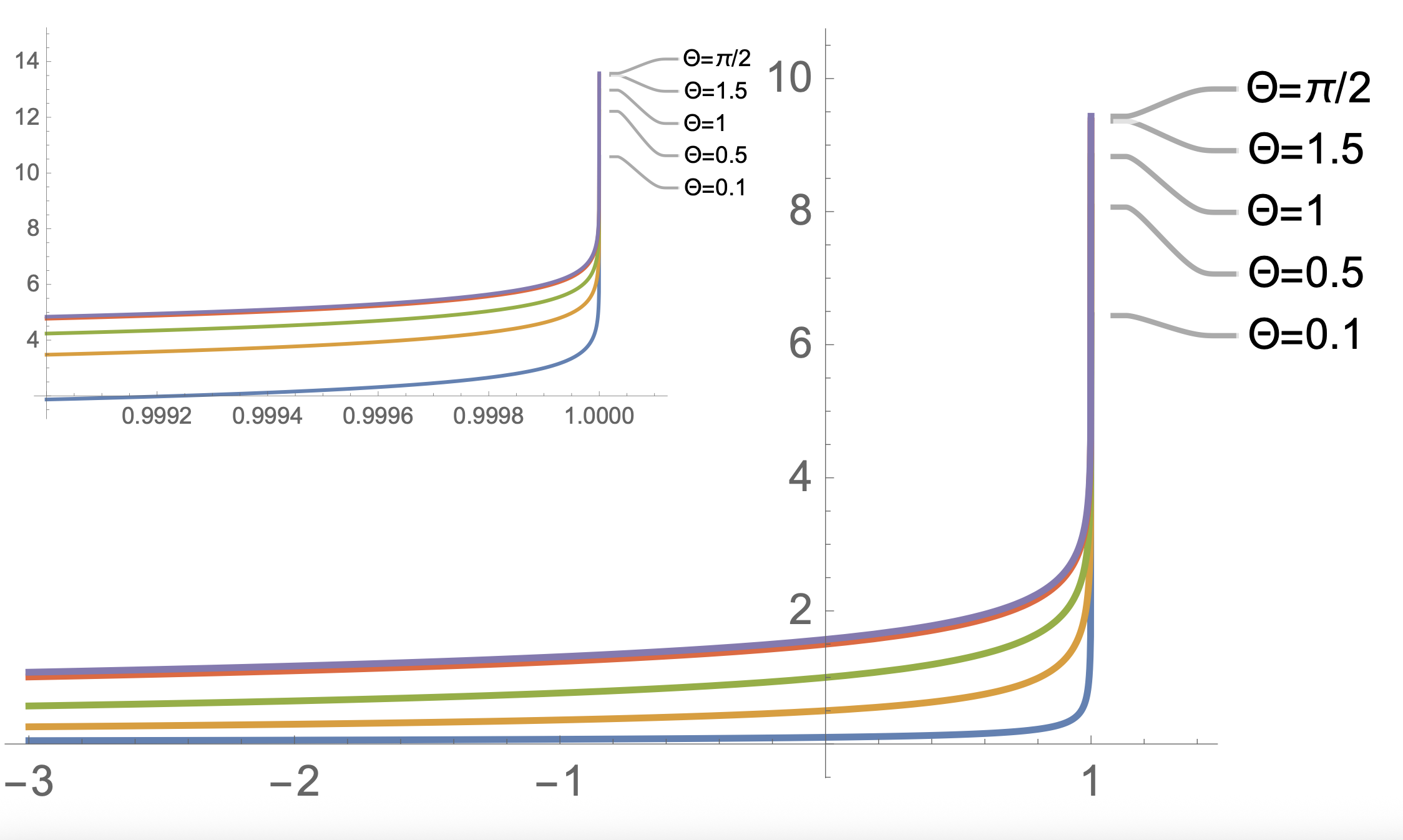}}\\
\end{center}
{\em The graph of $\frac{ F\left(\Theta\left|\frac{\alpha_\ell }{\alpha_\ell-1}\right.\right)}{\sqrt{1-\alpha_\ell } } $ in (\ref{alphaell}) for a range of $\Theta$. }

\medskip

The solution to (\ref{alphaell}) here is increasing in $\Theta$,  and increasing in $\alpha_\ell$ since the derivative $\frac{\partial}{\partial \alpha_\ell} \frac{ F\left(\Theta\left|\frac{\alpha_\ell }{\alpha_\ell-1}\right.\right)}{\sqrt{1-\alpha_\ell } }$ has the same sign as its numerator
\[ \alpha_\ell  \sin 2 \Theta  +\sqrt{2} ( \alpha_\ell -1) \sqrt{\frac{ \alpha_\ell  \cos  2 \Theta  + \alpha_\ell -2}{ \alpha_\ell -1}} \left(F\left[\Theta \left|\frac{\alpha }{  \alpha_\ell -1}\right.\right]-E\left[\Theta \left|\frac{ \alpha_\ell }{ \alpha_\ell-1}\right.\right]\right) \]

\begin{lemma}
\[ 2 \tan \Theta  \sqrt{1-\alpha_\ell  \cos ^2\Theta } - 2 \sqrt{1-\alpha_\ell } E\big[\Theta \left|\frac{\alpha_\ell }{\alpha_\ell -1} \right. \big]  = \int_{0}^{\Theta} \frac{  \tan ^2(\theta ) \;d\theta}{\sqrt{1-\alpha_\ell  \cos ^2(\theta )}} \]
\end{lemma}
\noindent{\bf Proof.} We calculate as follows.
\begin{eqnarray*}
\lefteqn{2 \tan \Theta  \sqrt{1-\alpha_\ell  \cos ^2\Theta } - 2 \sqrt{1-\alpha_\ell } E\big[\Theta \left|\frac{\alpha_\ell }{\alpha_\ell -1} \right. \big] } \\
&=&2 \tan \Theta  \sqrt{1-\alpha_\ell  \cos ^2\Theta } - 2  \int_{0}^{\Theta}  \sqrt{1-\alpha_\ell\cos^2(\theta) }   d\theta
\end{eqnarray*}
This difference is $0$ when $\Theta=0$. We differentiate this difference to obtain 
\[ \frac{2 \tan ^2(\Theta )}{\sqrt{1-\alpha  \cos ^2(\Theta )}} \geq 0 \] 
Hence we have
\begin{equation}
 \tan \Theta  \sqrt{1-\alpha_\ell  \cos ^2\Theta } -  \sqrt{1-\alpha_\ell } E\big[\Theta \left|\frac{\alpha_\ell }{\alpha_\ell -1} \right. \big]  = \int_{0}^{\Theta} \frac{  \tan ^2(\theta ) \;d\theta}{\sqrt{1-\alpha_\ell  \cos ^2(\theta )}}
\end{equation}
This establishes the lemma. \hfill $\Box$

\medskip

Next,  we consider the first term of (\ref{ineq}). We use (\ref{alphaell})
\begin{eqnarray*}
{2\ell}  \big(\sqrt{1-\alpha_\ell }+\frac{1}{\sqrt{1-\alpha_\ell }}\big) F\big[\Theta \left|\frac{\alpha_\ell }{\alpha_\ell -1}\right.\big]  
=\frac{2-\alpha_\ell}{2\pi}  \, \ell^2 \;m_\Omega  
\end{eqnarray*}
\begin{corollary}\label{cor} In Theorem \ref{thm2} we may replace the lower bound of (\ref{ineq}) by
\begin{eqnarray}\nonumber
 \int_A \IK(z,f) & \geq& \frac{2-\alpha_\ell}{4\pi}  \, \ell^2 \;m_\Omega +
 2\ell\int_{0}^{\Theta} \frac{  \tan ^2(\theta ) \;d\theta}{\sqrt{1-\alpha_\ell  \cos ^2(\theta )}} \end{eqnarray}
\end{corollary}

The two terms here have quite different behaviour as $\alpha_\ell$ varies in $(-\infty,1]$ and $\Theta$ varies in $[0,\frac{\pi}{2}]$.

\medskip

We now obviously need to understand the function $\alpha_\ell$.  However,  first we will identify the large round annuli in degenerating surfaces,  that is as $\ell\to0$.

\section{Collars in Riemann surfaces.}
 
 In this section we identify the annuli that we will need to investigate.  A {\em collar} (sometimes {\em cylinder}, of radius $\delta$ about a simple closed geodesic $\gamma$ in a Riemann surface $\Sigma$ is
 \[ A = \{z\in \Sigma: d_\Sigma(\gamma,z) < \delta\} \]
 Such a collar is isometric to the round annuli we considered above.  The existence of large,  that is $\delta$ large, collars about short geodesics is a well know property of hyperbolic manifolds in general \cite{GM} and yields the so called ``thick and thin'' or ``Margulis'' decomposition.  These are usually established by variants of J\o rgensen's inequality or related universal constraints. Here we want rather precise estimates and so we will use the following sharp universal constraint developed in Beardon's book \cite[\S 5.4]{Beardon}.
 
 \begin{lemma}\label{lemma1} Let $f$ and $g$ be hyperbolic M\"obius transformations whose axes cross at an angle $\theta$,  and that $\langle f,g\rangle$ is discrete and torsion free.  Then
\begin{equation} \sinh \frac{\tau_f}{2} \;  \sinh \frac{\tau_g}{2} \;  \sin \theta \geq 1. \end{equation}
 \end{lemma}
 The following corollary is what we seek.  It is known in various forms, see e.g. \cite{Buser}.
 \begin{lemma}\label{lemcol}
 Suppose that $\gamma$ is the shortest closed geodesic of length $\ell$ in a closed Riemann surface $\Sigma$.  Then $\gamma$ admits a collar of radius $\delta$, 
 \begin{equation}\label{delta}
 \delta = \sinh^{-1}\Big[\frac{1}{\sinh(\ell/2)}\Big]
 \end{equation}
 \end{lemma}
\noindent{\bf Proof.} We first note that an easy geometric argument shows the shortest geodesic $\alpha_\ell$ always has an embedded collar.  Let ${\mathcal C}_\delta$ be the largest such collar,  $\delta$ the radius, and $\beta$ the shortest geodesic passing across ${\mathcal C}_\delta$. Let $x$ be the hyperbolic length of $\beta$ and observe $x/2\geq \delta$. $\beta$ must intersect $\alpha_\ell$ at a point $x_0$ and suppose the angle of intersection is $\theta$. Lemma \ref{lemma1} gives 
 \[ \sinh \ell/2 \sinh x/2 \geq 1/\sin\theta \]
 Orthogonally project a point of $\beta\cap\partial{\mathcal C}_\delta$ (call it $x_1$) to $x_2\in\alpha_\ell$. Note the distance $d_{hyp}(x_0,x_2)\leq \ell/2$.  Some elementary hyperbolic trigonometry then gives (\ref{delta}). \hfill $\Box$
 
 \medskip
 
  We call the collar of Lemma \ref{lemcol} a {\em maximal} collar.  Note however it might not be the largest collar in a particular Riemann surface with a closed geodesic of length $\ell$.
  
We may use (\ref{area}) to establish the following.

\begin{corollary}\label{cor2} The maximal collar has hyperbolic area $2\ell/\sinh(\ell/2) \leq 4$,  and conformal modulus
\[ mod = \frac{ 4 \pi  {\rm sech}\left(\frac{\ell}{2}\right)}{\ell} \geq \frac{4 \pi }{\ell}-\frac{\pi \ell}{2} \]
\end{corollary}

\section{Bounds on $\alpha_\ell$} We seek good bounds as $\ell \to 0$. We consider $A$ as the maximal collar given to us from  Lemma \ref{lemcol} and which we use as data in Corollary \ref{cor}. Thus we must identify the behavour of the two terms there. 

We first make the following substitution
\begin{equation}\label{tdef} t = 1/\sqrt{1-\alpha_\ell} >0, \quad \alpha_\ell=1- t^{-2} ,\quad \frac{\alpha_\ell}{\alpha_\ell-1}=1-t^2\end{equation}
Next $t$ is defined by the relationship  (\ref{alphaell}) which now reads as
 \begin{equation}\label{trel}
 t F\left(\Theta\left|1-t^2\right.\right)  = \frac{\ell \;m_\Omega}{4\pi}
 \end{equation}
Then (\ref{trel}) defines a functional relationship between $t$ and $\ell$ as $ t F[\Theta;1-t^2]$ is concave and increasing with range $[0,\infty)$.  It is also immediate that $t\to0$ as $\ell\to0$. Then implicitly differentiating (\ref{trel}) w.r.t $\ell$, and setting $\ell=0$ gives 
\[ \frac{dt}{d\ell}\Big|_{\ell=0} = \frac{m_\Omega}{4\pi F[\Theta ;1]}\]
so that $t\ll\ell$ when $\Theta\approx\frac{\pi}{2}$,  and $t\gg\ell$ when $\Theta\approx 0$.

\medskip

Lemma \ref{lemcol} implies that for each $\ell>0$ we can find a round annulus $A$ in a surface with a core geodesic of length $\ell$ and radius $\delta$ satisfying the following.
\begin{eqnarray} \delta &=&  \sinh^{-1}\Big[\frac{1}{\sinh(\ell/2)}\Big]\approx \frac{2}{\ell}\label{smalldelta} \\
\Theta&=& \sin ^{-1}\Big(\frac{1}{\cosh \frac{\ell}{2} }\Big)=\frac{\pi }{2}-\frac{\ell}{2}+\frac{\ell^3}{48}-\frac{\ell^5}{768}+O\left(\ell^7\right)   \label{Thetaesti}
\end{eqnarray}

 \begin{lemma}\label{thetabounds} For $0< \ell\leq 1$, 
 \[ \frac{\pi }{2}-\frac{\ell}{2}\leq \Theta \leq  \frac{\pi }{2}-\frac{\ell}{2}+\frac{\ell^3}{48} \]
 \end{lemma}
\noindent{\bf Proof.} These follow as the series for $\Theta$ in $\ell$ is alternating,  the terms are decreasing and the series is uniformly convergent on $[0,1]$. \hfill $\Box$

\medskip
We suppose for each $\ell$ there is a mapping of finite distortion to another fixed (base) Riemann surface and
\[ f(A)=\Omega, \quad {\rm mod}(f(A))=m_\Omega \]
We have no a'priori information on  $m_\Omega>0$ as we do not know the (nontrivial) homotopy class in which $\Omega$ lies. However the base surface is not changing and so it has a maximal annular region of maximal modulus.  That is $m_\Omega$ is bounded above by a constant independent of $\ell$.  Also,  for each $\ell$ we may have a different mapping and image.  To recognise all this we set
\begin{equation}
m_\ell = \frac{1}{4\pi} m_\Omega.
\end{equation}
Now from (\ref{trel}) we have 
\begin{equation}\label{t2}\ell \,m_\ell = t F\left(\Theta\left|1-t^2\right.\right) = \int_{0}^{\Theta} \frac{t\,d\theta}{\sqrt{\cos^2\theta+t^2\sin^2\theta}}
\end{equation}
For a Riemann surface in a bounded part of moduli space,  the number $m_\ell$ is bounded above.  In particular on our reference surface $\Sigma_2$ we have $m_\ell$ is bounded above.
Now  the two terms we must examine from Corollary \ref{cor} are
\[ \frac{2-\alpha_\ell}{2\pi}  \, \ell^2 \;m_\Omega, \quad 
 4\ell\int_{0}^{\Theta} \frac{  \tan ^2(\theta ) \;d\theta}{\sqrt{1-\alpha_\ell  \cos ^2(\theta )}} \]
 \subsection{The second term.} With the data above we calculate
 First,  as $\ell m_\ell \to 0$ we have directly from (\ref{t2}) that $t\to 0$. So $\alpha_l\to -\infty$ and 
 \begin{eqnarray*}
  4\ell\int_{0}^{\Theta} \frac{  \tan ^2(\theta ) \;d\theta}{\sqrt{1-\alpha_\ell  \cos ^2(\theta )}} & \leq & \frac{4\ell}{\sqrt{1-\alpha_\ell}} (\tan(\Theta)-\Theta) \leq  \frac{4}{\sqrt{1-\alpha_\ell}} \to 0.
 \end{eqnarray*}
 Thus the second term will contribute little to the growth and so we ignore it.
  \subsection{The first term.}
  This term is $2(2-\alpha_\ell)  \, \ell^2 \;m_\ell$. Again, $m_\ell \ell \to 0$ and $\alpha_\ell \to -\infty$.  We need to show
  $ \frac{\ell^2}{t^2}\, m_\ell \to +\infty$ as $\ell\to 0$ given (\ref{t2}). Lemma \ref{thetabounds} and (\ref{t2}) together give
\begin{eqnarray}
  t F\left[\frac{\pi}{2} - \frac{\ell}{2} ;1-t^2 \right]      \leq    t F\left[\Theta;1-t^2 \right] \leq   t F\left[\frac{\pi}{2}- \frac{\ell}{2}+ \frac{\ell^3}{48} ;1-t^2 \right] \nonumber \\
\label{boundbyl}  t K[1-t^2] - \int_{0}^{\frac{\ell}{2}} \frac{t\,d\theta}{\sqrt{t^2 \cos^2\theta+\sin^2\theta}}       \leq   t F\left[\Theta;1-t^2 \right]   \leq  \\ \quad\quad t K[1-t^2]- \int_{0}^{\frac{\ell}{2}} \frac{t\,d\theta}{\sqrt{t^2 \cos^2\theta+\sin^2\theta}} + \int^{\frac{\ell}{2}}_{\frac{\ell}{2}-\frac{\ell^3}{48}} \frac{t\,d\theta}{\sqrt{t^2 \cos^2\theta+\sin^2\theta}}  
\end{eqnarray}
This now gives us a third order approximation in $\ell$.
\begin{eqnarray} 
   t K[1-t^2] -  F\left[\frac{\ell}{2};1-\frac{1}{t^2} \right]        \leq   t F\left[\Theta;1-t^2\right] \leq t K[1-t^2]- F\left[\frac{\ell}{2}  1-\frac{1}{t^2} \right]+  \frac{\ell^3}{48}  \end{eqnarray}
With the Elliptic $K$ function, $K[\cdot]=F[\frac{\pi}{2};\cdot]$.  The function $t\mapsto t K[1-t^2]$ is difficult to invert,  and so we seek a lower-bound on it obtained from a power series expansion
\[ t \log\frac{4}{t}+ \frac{1}{4}(-1 + \log\frac{4}{t}) t^3 \leq t K[ 1-t^2] \leq t\log\frac{4}{t}+ \frac{1}{2}(-1 + \log\frac{4}{t}) t^3 \] 
where we have made the assumption $t(-1 + \log\frac{4}{t})<\frac{1}{2}$ to get this.  This is implied by
\begin{equation}
\log\frac{4}{\ell} <1+\frac{1}{2\ell}, \quad \ell \leq -\frac{1}{2 W\left(-\frac{1}{8}\right)} = 0.15329\ldots
\end{equation}
Here $W$ is the Lambert function which we will consider more closely in a moment.  We have now achieved the following estimate.
\begin{lemma} For $0<t<\ell \leq -\frac{1}{2} W\left(-\frac{1}{8}\right)$ we have the estimate
\[   t \log\frac{4}{t} -\frac{\ell}{2} \leq t F\left(\frac{\pi}{2}-\frac{\ell}{2} \left|1-t^2\right.\right), \quad \ell \leq -\frac{1}{2 W\left(-\frac{1}{8}\right)}. \] 
\end{lemma}
Our estimates above are finer than this,  but do not measurably improve the result we obtain below.
\begin{center}
\scalebox{0.33}{\includegraphics{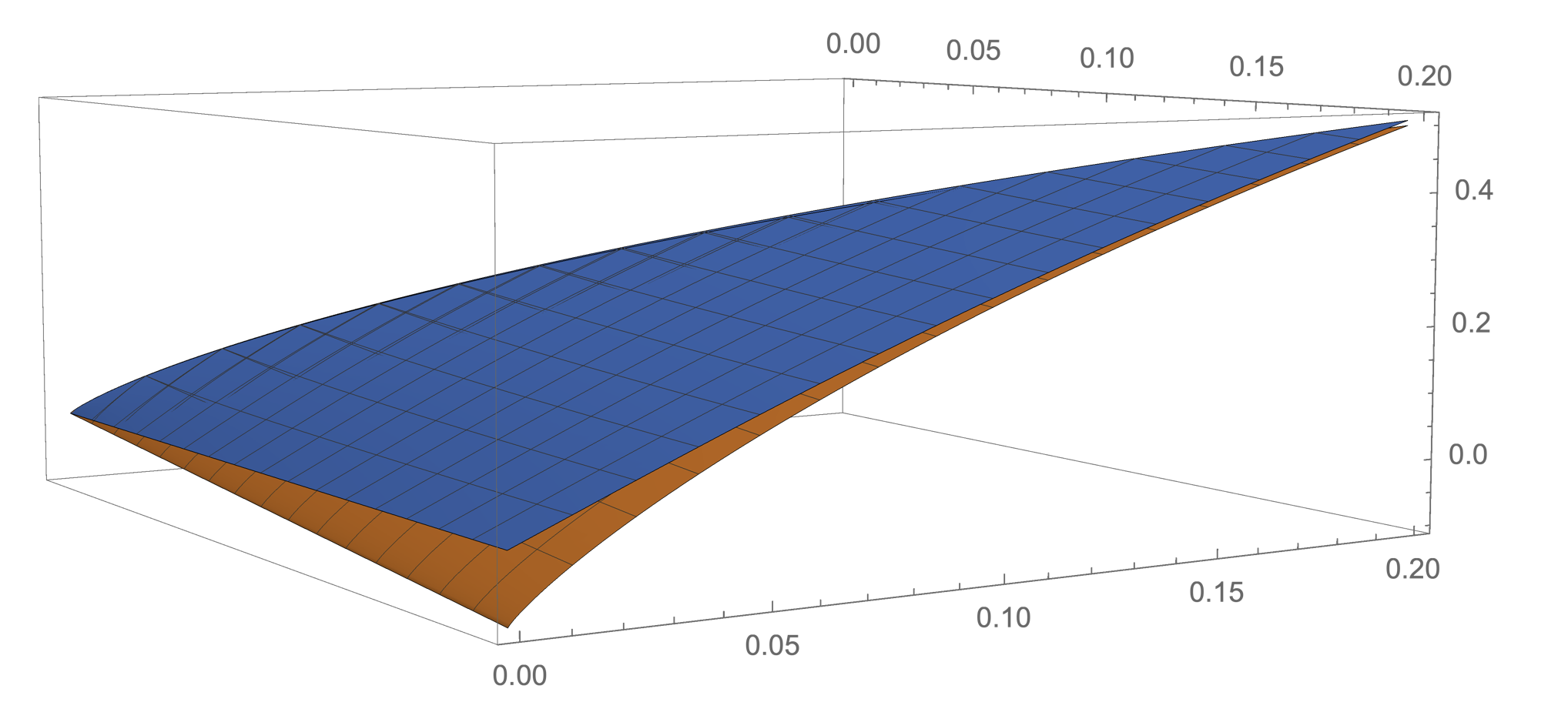}}
\end{center}
\noindent{\bf Figure X} {\em The graphs of $t F\left(\frac{\pi}{2}-\frac{\ell}{2} \left|1-t^2\right.\right)$ (top) and  $t \log\frac{4}{t} -\frac{\ell}{2}$ (below).}

Now given $\ell\in [0,0.15]$ we choose $t_0$ so that
\begin{equation}\label{tochoose}
\ell\,m_\ell= t_0\log\big[\frac{4}{t_0}\big]-\frac{\ell}{2}, \quad t_0=4 e^{W_{-1}\big[ -\frac{1}{4}\, \ell (m_{\ell}+\frac{1}{2})\big] }\end{equation}
Here $W$ is the Lambert $W$,  or ProductLog function,  the principal solution of $z=we^w$ in $\mathbb C$.  Actually we have used $W_{-1}$,  the other real branch. Remarkably good asymptotics for this function are given in \cite{BLJ} and explicit bounds in \cite{chat}.  From these we have the following when $ \ell m_{\ell}\geq   4 e^{-1}$.
\begin{eqnarray*}
t&=&4 e^{W_{-1}\big[ -\frac{1}{4}\, \ell (m_{\ell}+\frac{1}{2})\big] },\\
 &\approx& \frac{\ell (m_\ell+\frac{1}{2})}{ e^{2/\beta }}    \exp \left(\frac{2 \sqrt{2}}{\sqrt{2} \beta +\beta ^2 \sqrt{\log \left(\frac{4} {(m_{\ell }+\frac{1}{2})\, \ell }\right)-1}}\right), \quad \beta=0.3205.
\end{eqnarray*}
 This approximation has maximum relative error of $3.8\times 10^{-3}$.  However these approximations clearly do not have the property that we have already identified as necessary for us,  namely that $t/\ell \to 0$ as $\ell\to 0$. 
 
 \medskip
 
 Since $t F\left(\frac{\pi}{2}-\frac{\ell}{2} \left|1-t^2\right.\right)$ is decreasing we have the following theorem.
 \begin{theorem} For $0<\ell<0.15$, let $\Theta=\sin^{-1}(1/\cosh(\ell/2))\approx \frac{\pi}{2}-\frac{\ell}{2}$ and $t$ be the solution to the equation
 \[t F\left(\Theta\left|1-t^2\right.\right) = \ell m_\ell\]
 Then 
 \[ t<t_0=4 e^{W_{-1}\big[ -\frac{1}{4}\, \ell (m_{\ell}+\frac{1}{2})\big] }\]
 \end{theorem}
 
 Now we consider the term $2(2-\alpha_\ell) \, \ell^2 \;m_\ell$ that we must bound from below.
\begin{eqnarray*} 
2(2-\alpha_\ell) \, \ell^2 \;m_\ell & = & 2(1+\frac{1}{t^2})  \, \ell^2 \;m_\ell = 2 t \ell (1+\frac{1}{t^2})  \, \frac{\ell \;m_\ell}{t}\\
& = & 2 \ell (t+\frac{1}{t})  \,  F\left(\Theta\left|1-t^2\right.\right) \quad \mbox{using (\ref{t2})} \\
& \geq & 2 \ell (t_0+\frac{1}{t_0})  \,  F\left(\Theta\left|1-t_0^2\right.\right) \geq   \frac{2\ell}{t_0}   \,  F\left(\Theta\left|1-t_0^2\right.\right) \\
&\geq & \frac{2\ell}{t_0}   \,  F\left(\frac{\pi}{2}-\frac{\ell}{2}\Big|1-t_0^2 \right) \\
\end{eqnarray*}
Now an application of L'H\^{o}spital's rule gives
\[ \lim_{\ell\to 0} \frac{t_0}{\ell} =\lim_{\ell\to 0} \frac{-2 m-1}{2 \left(W_{-1}\left(-\frac{1}{8} l (2 m+1)\right)+1\right)} = 0 \]
 and as 
 \[ F\left(\frac{\pi}{2}-\frac{\ell}{2}\Big|1-t_0^2 \right) = \int_{0}^{\frac{\pi}{2}-\frac{\ell}{2}} \frac{d\theta}{\sqrt{\cos^2(\theta)+t_0^2\sin^2(\theta) }} \to \infty, \quad  \ell \to 0 \]
 we see that  
 \begin{equation}\label{fterm}
 2(2-\alpha_\ell) \, \ell^2 \;m_\ell  \to +\infty, \quad \ell \to 0.
 \end{equation}
 This is all we need for properness in this direction but we would like to identify explicit lower bounds so as to compare with the Teichm\"uller distance.
 
 \begin{theorem}\label{thm7}  Let $\Sigma_1, \Sigma_2$ be Riemann surfaces of genus $g$. Let
$m = \frac{1}{2\pi} \max_A   mod(A)$ where $A$ is an annular region about a simple closed geodesic of $\Sigma_2$.
Suppose that $\Sigma_1$ has a shortest geodesic of length 
\[  \ell \leq \frac{8}{(2m+1)e^{2(2m+1)} }\] 
If $f:\Sigma_1\to \Sigma_2$ is a homeomorphism of finite distortion,  then
\begin{equation}\label{ff}
\int_{\Sigma_1} \IK(z,f) \; d\sigma(z) \geq    4\pi(g-1)-4+ \frac{\sqrt{2}}{2m+1}\Big(1+\log\frac{4}{\ell}\Big){  \log  \frac{8}{\ell (2m+1)}}    \end{equation}
\end{theorem}
\noindent{\bf Proof.} We must make a  first observation.  In order to achieve a ``clean'' formula for $t_0$ we bounded the term 
\[  \int_{0}^{\frac{\ell}{2}} \frac{t\,d\theta}{\sqrt{t^2 \cos^2\theta+\sin^2\theta}}  = t F[\frac{\ell}{2};1-t^2]  < \frac{\ell}{2} \]
 This was reasonable in that situation,  but here we now know $t\ll \ell$. This error term introduced the term $m+\frac{1}{2}$ instead of just $m$ in our formula for $t_0$.  This is of no real consequence if the number $m$ is bounded below. However if $m\to 0$,  which might occur if the image $f(A)$ is concentrating on a homotopically nontrivial simple closed curve, then the bounds on $t_0$ are simply not good enough to use directly. 
\begin{center}
\scalebox{0.25}{\includegraphics{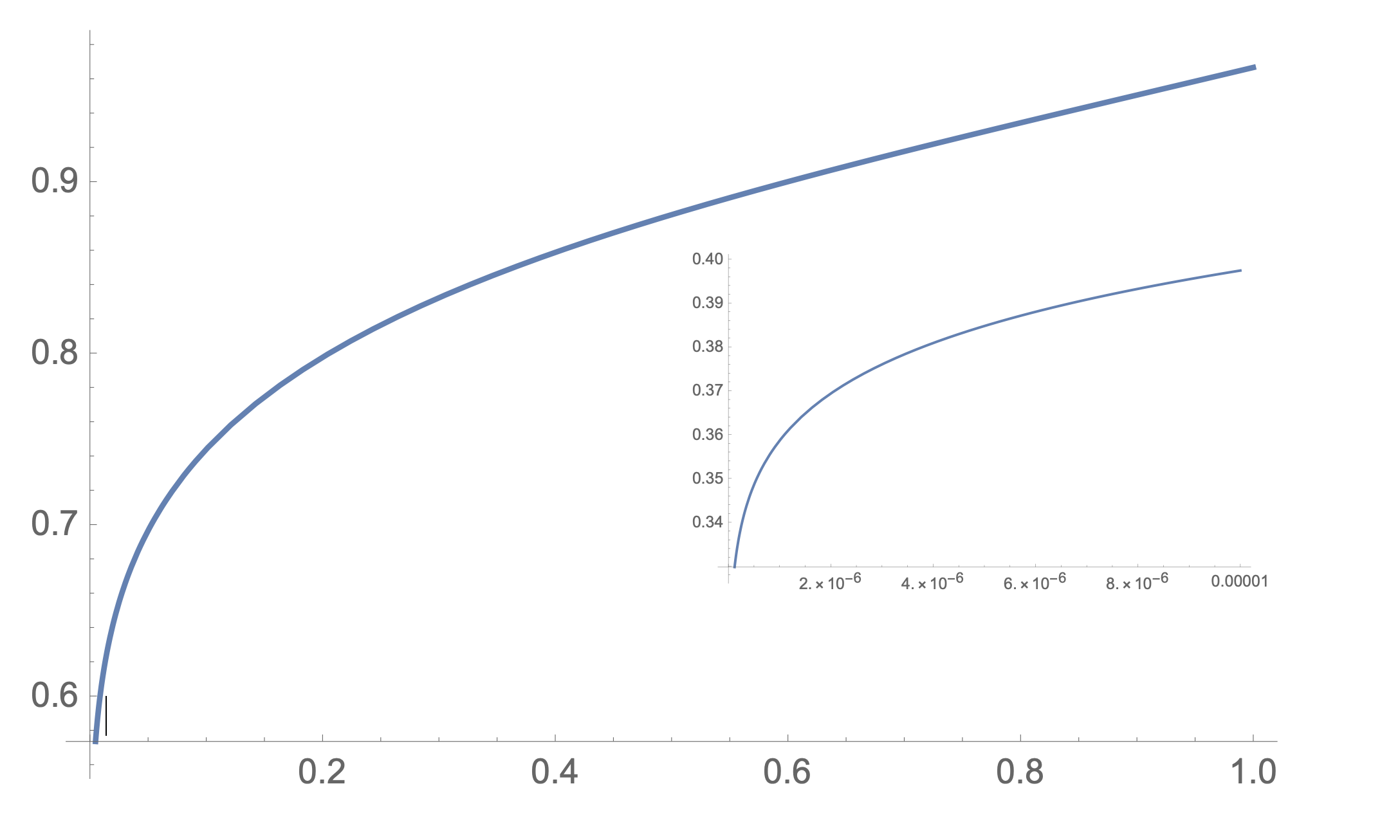}}
\end{center}
\noindent {\bf Figure 1} {\em The ratio of the two bounds when $t\approx \ell/\log(4/\ell)$.  The graph of the ratio $\frac{2}{\ell} F\big(\frac{l}{2}\big|1-{l^{-2}}\log ^2\left(\frac{4}{l}\right)\big) $, inset with $\ell \in [0,10^{-5}]$.  We used the bound $ratio \leq 1$}.

\medskip
Thus, simply putting our bound for $t_0$ in  $ \frac{\ell^2}{t_0^2}\, m_\ell$ will not give us divergence to $+\infty$ as $\ell\to0$, if $m\to 0$ also. Recall that $t_0$ is defined in terms of $m$. The way forward is simply to follow (\ref{fterm}).
 \begin{eqnarray*}
 2(2-\alpha_\ell) \, \ell^2 \;m_\ell  & = & 2 \ell (t+\frac{1}{t})  \,  F\left(\Theta\left|1-t^2\right.\right) \quad \mbox{using (\ref{t2})} \\
& \geq &2 \ell (t_0+\frac{1}{t_0})  \,  F\left(\frac{\pi}{2}-\frac{\ell}{2}\Big|1-t_0^2 \right) 
  \geq   \frac{2\ell}{t_0}   \,  F\left(\frac{\pi}{2}-\frac{\ell}{2}\Big|1-t_0^2 \right)
 \end{eqnarray*}

\begin{lemma}\label{t0est} For $x\in [0,8/e]$ we have the inequality.
\begin{equation}\label{ineq2} 4 \exp\big(W_{-1}\big[ -\frac{x}{8}\big]\big) \leq \frac{x}{ 2\log  \frac{8}{x}}
\end{equation}
\end{lemma}
\noindent{\bf Proof.}  Both functions vanish at $0$, and agree at $x=1/e$,  the singularity of the Lambert function. 
The derivative of the difference
\[ \frac{d}{dx} \left(4 \exp\big(W_{-1}\big[ -\frac{x}{8}\big]\big) - \frac{x}{ 2\log  \frac{8}{x}}\right) = \frac{1}{2} \left(\frac{\log \left(\frac{x}{8}\right)-1}{\log ^2\left(\frac{8}{x}\right)}-\frac{1}{W_{-1}\left(-\frac{x}{8}\right)+1}\right) \]
is initially negative so the difference between these functions is initially increasing. The derivative of the left-hand side of (\ref{ineq2}) is negative and decreasing while the derivative of the right-hand side is positive and increasing. Thus there is at most one critical point and the result follows. \hfill $\Box$

\medskip The graphs of these two functions and there difference is illustrated below in Figure 2

\begin{center}
\scalebox{0.35}{\includegraphics{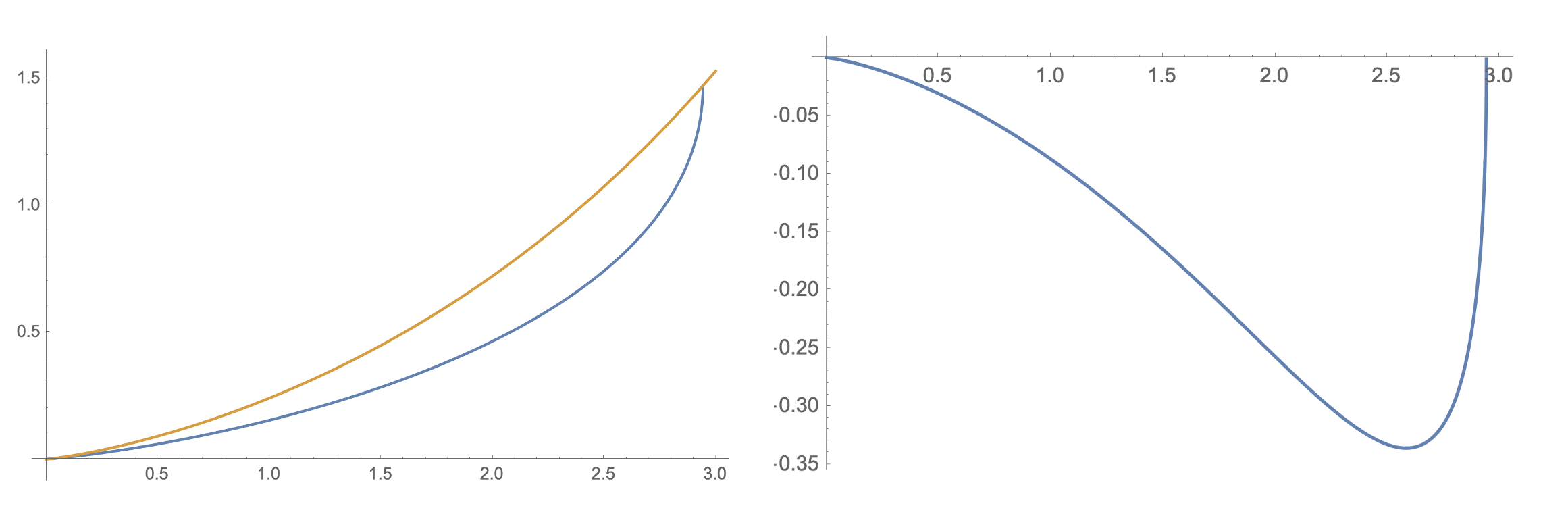}}
\end{center}
\noindent{\bf Figure 2} The graphs of the functions and their difference for Lemma \ref{t0est}.

\medskip

Now Lemma \ref{t0est} gives us that for $\ell \leq \frac{1}{e (2m+1)}$,
\begin{equation}
\frac{\ell}{t_0} \geq \frac{2}{2m+1}{  \log  \frac{8}{\ell (2m+1)}}
\end{equation}
Recall here that $m$ is data from the image and is uniformly bounded above.

The next term we must consider is   $ F\left(\frac{\pi}{2}-\frac{\ell}{2}\Big|1-t_0^2 \right)$.  We have from Lemma \ref{t0est}
\begin{eqnarray*}
F\left(\Theta\Big|1-t_0^2 \right) & = & \int_{0}^{\Theta} \frac{d\theta}{\sqrt{1-(1-t_0^2)\sin^2\theta}} = \int_{0}^{\Theta} \frac{d\theta}{\sqrt{\cos^2\theta+t_0^2\sin^2\theta}} \\
& = & \int_{0}^{\pi/4} \frac{d\theta}{\sqrt{\cos^2\theta+t_0^2\sin^2\theta}} +\int^{\Theta}_{\pi/4} \frac{d\theta}{\sqrt{\cos^2\theta+t_0^2\sin^2\theta}} 
\end{eqnarray*}
This term is strongly asymptotic to $\log\frac{4}{\ell}$, as a power series analysis shows,  but getting explicit bounds seems rather challenging. We use the following lemma to continue this calculation. 

\begin{lemma} For $0<\ell < \frac{8}{(2m+1)e^{2(2m+1)} }$ and $\theta\in [\frac{\pi}{4},\frac{\pi}{2}-\frac{\ell}{2}]$
\[ \cot(\theta) \geq t_0 \]
\end{lemma}
\noindent{\bf Proof.} The term $\cot(\theta)$ is strictly decreasing and so we need only prove 
\[  \cot \big(\frac{\pi}{2}-\frac{\ell}{2}\big) = \tan \frac{\ell}{2}  \geq \frac{\ell(2m+1)}{\log  \frac{8}{\ell (2m+1)}} \geq t_0 \]  
The latter inequality is Lemma \ref{t0est}. As $\tan(\frac{\ell}{2}) \geq \frac{\ell}{2}$ it is now enough to ensure
\[  \frac{1}{2}  \geq \frac{2m+1}{\log  \frac{8}{\ell (2m+1)}}  \]
which was our choice for the bound on $\ell$. \hfill $\Box$

\medskip This result implies that $t_0^2\sin^2\theta \leq \cos^2 \theta$ in the domain of the following integral. Now we calculate
\begin{eqnarray*}
F\left(\Theta\Big|1-t_0^2 \right) &\geq  & \frac{\sqrt{2}}{\sqrt{t_0^2+1}}+\int^{\Theta}_{\pi/4} \frac{d\theta}{\sqrt{\cos^2\theta+t_0^2\sin^2\theta}} \\
&\geq  & \frac{\sqrt{2}}{\sqrt{t_0^2+1}}+\frac{1}{\sqrt{2}} \int^{\frac{\pi}{2}-\frac{\ell}{2}}_{\pi/4} \frac{d\theta}{\cos\theta} \\
&\geq  & \frac{\sqrt{2}}{\sqrt{t_0^2+1}}+ \log  \cot \big(\frac{\ell}{4}\big) -\sinh ^{-1}(1) \\
&\geq  &  \frac{1}{\sqrt{2}}\Big(1+\log\frac{4}{\ell}\Big).
\end{eqnarray*}

 This final lower bound, together with Lemma \ref{t0est}, completes the proof of Theorem \ref{thm7}.
  \hfill $\Box$
 
 \subsection{An upper bound.}
Now we seek an upper bound. If $m_\Omega$ is the fattest annulus in the target, then there is a quasiconformal map $f:\IA\to f(\IA)$,  $mod(f(\IA))\leq m_\Omega$ of distortion 
\[ \IK(z,f_0) \leq \frac{1}{2}\Big(\frac{\mod(\IA_1)}{m_\Omega} +\frac{m_\Omega}{\mod(\IA_1)}\Big) \]
and so for the extremal in $L^1$ we must have
\begin{eqnarray*}
\inf_f \int_{\IA_1} \IK(z,f) d\sigma  & \leq & \inf \int_{\IA_1} \IK(z,f_0) d\sigma = \IK(z,f_0) |\IA_1|_{hyp} \\
& \leq & \frac{1}{2}\Big(\frac{\mod(\IA_1)}{m_\Omega} +\frac{m_\Omega}{\mod(\IA_1)}\Big)\frac{ 2\ell}{\sinh(\ell/2)} 
\end{eqnarray*}
 This actually gives a lower bound on the Teichm\"uller distance and Kerchoff's extensive investigation into the relationship between extremal modulus and the Teichm\"uller distance \cite{Kerchoff} shows this to be of the correct order.
 
 Choosing $\ell$ small enough so that ${mod(\IA_1)}>m_\Omega$ for simplification gives
 \[ \inf_f \int_{\IA_1} \IK(z,f) d\sigma \leq \frac{4 \pi}{\sinh(\ell)},\quad m_\Omega<mod(\IA_1)\]

  \section{Proof for Theorem \ref{thm1}{\bf (b)} }
  
  In this case $\Sigma_1$ is bounded in Teichm\"uller space and so the hyperbolic metric of $\Sigma_1$ is uniformly bounded above and below. Thus we cannot exploit any degeneration in the metric to provide bounds, and the case becomes more or less the same as the Euclidean case.  We will in essence  reduce to this case.
  
  \medskip
  
  We begin with the following lemma.
 \begin{lemma}\label{lem1}
 Let $(\Sigma,\sigma)$ be an analytically finite Riemann surface with hyperbolic metric $\sigma$ and $\Omega$ an annular region lying in the homotopy class of $\gamma$, a simple closed geodesic of length $\ell>0$ (which need not meet $\Omega$). Then $\Omega$ is isometric to some homotopically nontrivial annular region $\tilde{\Omega} \subset A(1/s,s)$ with $s=e^{\pi^2/\ell}$. Let
 \[ \lambda_0 = \sup_{z\in \Omega} \rho_\Sigma(z) \]
 where $\rho_\Sigma$ is the hyperbolic metric density of $\Sigma$. Then
\begin{equation}\label{2a} 
\tilde{\Omega} \subset A(1/r,r), \quad {\rm where} \quad\lambda_0 = \frac{\ell}{2\pi}\; \frac{1}{r \cos\Big(\frac{\ell \log r }{2\pi}\Big)}, \quad  1<r\leq   s=e^{\pi^2/\ell}.
\end{equation}
 \end{lemma}
 \noindent{\bf Proof} The first part is simply the uniqueness statement of the Riemann mapping theorem for doubly connected domains.  We now use (\ref{4a}) to solve 
   \begin{equation} 
\lambda_0 = \frac{\ell}{2\pi}\; \frac{1}{r \cos\Big(\frac{\ell \log r }{2\pi}\Big)}
\end{equation}
for $r$. The function on the right-hand side is strictly increasing  in $r$  to $+\infty$ for $1<|z|<s=e^{\pi^2/\ell}$.  That $\tilde{\Omega}$ is isometrically embedded tells us that the maximum value of the hyperbolic density is also at most $\lambda_0$ and we deduce (\ref{2a}).  \hfill $\Box$

\medskip
Now unwinding definitions enables us to consider obtaining estimates in the planar case.

\begin{corollary}\label{cor3} For any function $u:\Omega \to \IR$, and with the obvious notation, 
\[ \int_\Omega  u(z) d\sigma(z) = \int_{\hat{\Omega}}  \hat{u}(z) \; \lambda_{A_s}(z) dz \leq \lambda_0 \int_{\hat{\Omega}}  \hat{u}(z) \]
\end{corollary}

The previous lemma and its corollary allow us to normalise the situation we find ourselves in.  This we now discuss.

\medskip

We have a surface $\Sigma$ of genus $g\geq 2$ and hyperbolic metrics $\Sigma_1=(\Sigma,\sigma_1)$ and $\Sigma_\ell=(\Sigma,\sigma_\ell)$,   a homeomorphism $f_\ell:\Sigma\to \Sigma$,  $\ell\in [0,1]$,  of finite distortion between them in a fixed homotopy class $[f]$ of homeomorphism. There is a simple closed geodesic $\alpha_\ell \subset \Sigma_\ell$ of length at most $\ell$,  and for each $\ell$, $f_\ell^{-1}(\alpha_\ell)$ is in the homotopy class of a fixed simple closed geodesic  $\alpha\in \Sigma_1$. Since typically $[f]\neq[Identity]$ $\alpha$ and $\alpha_\ell$ might be quite ``different'' geodesics.   The geodesic $\alpha_\ell$ is the central core of a fat annulus $A_\ell$ (the collar of Corollary \ref{cor2}),  and 
\[ mod(A_\ell) \geq  \frac{ 4 \pi  {\rm sech}\left(\frac{\ell}{2}\right)}{\ell} \geq \frac{4 \pi }{\ell}-\frac{\pi \ell}{2}  \]

The geometry of $\Sigma_1$ is unchanging as $\ell$ varies,  and $\alpha$ also has a maximal  (round)  collar,  say $\Omega$ - this is fixed and independent of $\ell$. The following lemma gives explicit sharp bounds in the nicest situation.

\begin{lemma} Suppose $f^{-1}(A_\ell) \subset \Omega$.  Let $\lambda_0 = \min_{\Sigma_1} |\sigma_1|$. Then
\begin{equation}\label{mest1}
\int_{\Sigma_1} \IK(z,f)\; d\sigma_1 \geq \lambda_0 |\Omega| \left( \frac{mod(A_\ell)}{mod(\Omega)} +  \frac{mod(\Omega)} {mod(A_\ell)}\right) 
\end{equation}
\end{lemma}
\noindent{\bf Proof.} Since $\Omega$ is a round annulus, and $A_\ell\subset f(\Omega)$ we have $mod(f(\Omega))\geq mod(A_\ell)$.  Let $\lambda = \min_\Omega \sigma_1$, note $\sigma_1$ is a conformal metric so this is well defined and $\lambda>0$, Then with the obvious notation,
\begin{eqnarray*} \int_{\Sigma_1} \IK(z,f)\; d\sigma_1& \geq & \int_{\Omega} \IK(z,f) \; d\sigma_1 \geq \lambda \int_{\hat{\Omega}} \IK(z,\hat{f}) \; dz, \end{eqnarray*}
where we have used the isometric embedding provided by Corollary \ref{cor3}. The right-hand side of (\ref{mest1})  is increasing in $mod(A_\ell)$,  so the result follows from \cite[Theorem 2]{AIM2}. \hfill $\Box$

\medskip
Note that the right-hand side here of (\ref{mest1})  is  
$\approx   \frac{\lambda_0 |\Omega| }{mod(\Omega)} \frac{1}{\ell}$.

\medskip

In fact the result \cite[Theorem 2]{AIM2} only requires that $f^{-1}(A_\ell)$ lies in some round annulus $\Omega$,  it need not have anything to do with $\alpha$ - but that it is a round annulus is key.  Our issue in general is that this might never be the case,  $\alpha$ could be a very long simple closed geodesic without any reasonable collar in its homotopy class. However it is fixed and so we set 
\[ r = \max\{t: \mbox{ $\alpha$ lies in a round annulus of width $t$} \}\]
We do not ask that $\alpha$ lies in the middle  of this annulus.

We proceed as follows.  First, $f^{-1}(A_\ell)$ is in the same homotopy class,  hence the same isotopy class,  as $\Omega$ and thus there is an annular region $\hat{\Omega}$ containing $\Omega$ and $f^{-1}(A_\ell)$.  Then if $\beta$ is a path between the two components of $\hat{\Omega}$ its length is at least $2r$.  With $\lambda_+$ and $\lambda_-$ being the maximum and minimum of the conformal density of the hyperbolic metric on $\hat{\Omega}$, and after identifying these things with there isometric imaged in annuli in $\IC$ as per Lemma \ref{lem1} we find ourselves in the situation that
\begin{itemize}
\item $\hat{\Omega}$ is an annular region in $\IC$,
\item any path between the two boundary components of $\hat{\Omega}$ has (Euclidean) length at least  $\frac{2r}{\lambda_+}$.
\end{itemize}
We now turn to consider standard estimates on the moduli of path families - a conformal invariant based on the length-area methods. Thus
\begin{itemize}
\item $\frac{\lambda_+}{2r}$ is admissible to the family of curves $\Gamma$ joining the components of $\partial\hat{\Omega}$.  That is 
\[ \int_\gamma \frac{\lambda_+}{2r} |dz| = \frac{\lambda_+}{2r}\int_\gamma  |dz| = length(\gamma)\frac{\lambda_+}{2r} \geq 1\]
\item $mod(\hat{\Omega})\leq \int_{\hat{\Omega}} \Big(\frac{\lambda_+}{2r}\Big)^2 |dz|^2$, as the modulus is the infimum of the integral of the squares of admissible functions.
\item By change of variables,
\[ \int_{f(\gamma)} \|Dh\| \frac{\lambda_+}{2r}  |dw| \geq 1\]
that is $\|Dh\|   \frac{\lambda_+}{2r} $ is admissible to the curve family $f(\Gamma)$.
\item Hence 
\begin{eqnarray*} \int_{ f(\hat{\Omega})} \|Dh\|^2 \Big(\frac{\lambda_+}{2r}\Big)^2 \geq mod(f(\hat{\Omega}))\geq mod(A). 
\end{eqnarray*}
as $A\subset f(\hat{\Omega})$.  
\item The change of variables formula gives
\[ \int_{ f(\hat{\Omega})} \|Dh(w)\|^2 |dw|^2 = \int_{\hat{\Omega}} \IK(z,f)\; |dw|^2   \]
and hence
\[ \int_{\hat{\Omega}} \IK(z,f)\; \lambda_-^2 \, |dw|^2 \geq   \Big(\frac{2r\lambda_-}{\lambda_+}\Big)^2 mod(A) \]
\item Hence
\[ \int_{\Sigma_1} \IK(z,f)\; d\sigma \geq \Big(\frac{2r\lambda_-}{\lambda_+}\Big)^2 mod(A) \] 
\end{itemize}

\begin{theorem} Let $U\subset \IC$ be a ring for which every homotopically nontrivial curve has length $d>0$. Let $\tilde{f}:U \to V$ be a homeomorphism of finite distortion. Then we have the following upper bound on the conformal modulus of $V$.
\begin{equation}
\frac{1}{d^2} \int_U \IK(z,f) \geq mod(V).
\end{equation}
\end{theorem}

\section{The $L^p$ case.}

As with the calculations above in the case that $p=1$, the case that the domain is bounded in moduli space offers little new to consider as the hyperbolic metric is bounded above and below and there are no large annuli in $\Sigma_1$.  Thus this case is basically the same, albeit a little more complicated, as the Euclidean case and so we do not pursue it here.  We consider the case that presented the most difficulty above,  when $\Sigma_1$ had a very short geodesic of length $\ell$ and $\Sigma_2$ is fixed and of bounded geometry. Thus we proceed as above and identify a large round annulus $\IA_1$ in $\Sigma_1$ of modulus at least $\frac{4 \pi }{\ell}-\frac{\pi \ell}{2}$ by Corollary \ref{cor2}. We may assumme the target is a round annulus of modulus $m$ bounded above (and possibly varying with $\ell$).
 
The equation (\ref{ve}) with $\Psi=t^{p}$, $p>0$, reads as 
 
\begin{equation}\label{uxeqn2}
\Big(1-\frac{1}{u_x^2}\Big)(u_x+\frac{1}{u_x})^{p-1} = \frac{\alpha}{\lambda(x)}  .
\end{equation}
where we have absorbed the constant $p$ into $\alpha$.  In our situation $\Sigma_1$ has a relatively much shorter geodesic, with a much larger collar and hence $u_x\ll 1$.  We then consider the   approximating equation (\ref{uxeqn}).
\begin{equation} \label{approxuxeqn}
  u_x \ll  1, \quad  u_x^{-(p+1)} = \frac{\alpha}{\lambda(x)}, \quad u_x= \Big(\frac{\lambda(x)}{\alpha}\Big)^\frac{1}{p+1} .
\end{equation}
We discuss how good this approximation is later.  As $u_x>0$ we also absorb the $-1$ term into $\alpha$,  so $\alpha>0$.
We recall (\ref{realweight}),
\[ 
\lambda(x) = \frac{\ell^2}{\cos^2\Big(\ell (x -\frac{a}{2})\Big)}, \quad  0\leq x \leq a, \]
where for simplicity of notation we have set $a=\frac{1}{2\pi}mod(\IA_1)\gg1$ and $b=\frac{1}{2\pi}mod(\IA_2)\ll a$, bounded. 

Next,  we have the boundary conditions determining $\alpha$. We see that $2\pi b$ is modulus of the image $f(\IA_1)$.  This is functionally dependent on $\ell$ but is bounded above.
\begin{eqnarray*}
2\pi b&=&\int_0^au_xdx=\int_0^a\left(\frac{\lambda(x)}{\alpha}\right)^{\frac{1}{p+1}}dx\\
\alpha&=&\frac{1}{(\pi b)^{p+1}}\left(\int_0^{a/2}\lambda(x)^{\frac{1}{p+1}}dx\right)^{p+1}\\
&=&\frac{1}{(\pi b)^{p+1}}\ell^{1-p}\left(\int_0^{a\ell/2}\frac{dx}{\cos^{2/(p+1)}(x)}\right)^{p+1}\\
\alpha^{p/(p+1)}&=&\frac{1}{(\pi b)^p}\ell^{(1-p)p/(p+1)}\left(\int_0^{a\ell/2}\frac{dx}{\cos^{2/(p+1)}(x)}\right)^p
\end{eqnarray*}

As earlier,  we calculate the integral of distortion as follows.
\begin{eqnarray*}  
\int_{\IA_1} \IK^{p}(z,f) d\sigma_1  & \geq &   \int_{0}^{a} u_x^{-p} \lambda(x) dx \\
& = &   \int_{0}^{a} \Big(\frac{\alpha}{\lambda(x)}\Big)^{p/(p+1)} \lambda(x) dx = {\alpha}^{p/(p+1)} \int_{0}^{a}  \lambda(x)^{1/(p+1)} dx \\ 
& = &\ell^{2/(p+1)} {\alpha}^{p/(p+1)} \int_{0}^{a}  \frac{dx}{\cos^{2/(p+1)}\Big(\ell (x -\frac{a}{2})\Big)}  \\
 \\
 & = &2\ell^{1-p} \frac{1}{( \pi b)^p} \Big(\int_{0}^{a\ell/2} \frac{dx}{\cos^{2/(p+1)} (x)}\Big)^{p+1}    \\
 \\
&=& \frac{ 2\ell^{1-p}}{( \pi b)^p}  \Big[-\frac{1}{2} B_{\cos ^2(x)}\left(\frac{1}{2}-\frac{1}{p+1},\frac{1}{2}\right)\Big] \Big|_{x=0}^{x=a\ell/2}
\end{eqnarray*}
In terms of the Euler $\beta$-function.   We have already identified an estimate for $a\ell$ at Corollary \ref{cor2}. 
\begin{equation}\label{al}
a  = \frac{1}{4\pi}mod(\IA_1) = \frac{1}{4\pi} \frac{4\pi}{\ell} \sin^{-1}(\tanh \delta) \approx  \frac{1}{\ell} \Big(\frac{\pi}{2}-\frac{\ell}{2}\Big),
\end{equation}
and hence $a\approx \frac{\pi}{2\ell}-\frac{1}{2}$ by (\ref{Thetaesti}). 

Then 
\begin{eqnarray*} \lefteqn{ \Big[-\frac{1}{2} B_{\cos ^2(x)}\left(\frac{1}{2}-\frac{1}{p+1},\frac{1}{2}\right)\Big] \Big|_{x=0}^{x=a\ell/2} }  \\
& = &\frac{1}{2} \left(B_0\left(\frac{1}{2}-\frac{1}{p+1},\frac{1}{2}\right)-B_{\sin ^2\left(\frac{l+\pi }{4}\right)}\left(\frac{1}{2}-\frac{1}{p+1},\frac{1}{2}\right)\right)\\ 
&\approx&\frac{1}{2} \left(B_1\left(\frac{1}{2}-\frac{1}{p+1},\frac{1}{2}\right)-B_{\frac{1}{2}}\left(\frac{1}{2}-\frac{1}{p+1},\frac{1}{2}\right)\right)-\ell 2^{\frac{1}{p+1}-2}+\frac{\ell^2 2^{\frac{1}{p+1}-4}}{p+1} +O\left(\ell^3\right)
\end{eqnarray*}
The constant term here is increasing in $p$ and 
\[  \frac{\pi}{4}  \leq \frac{1}{2} \left(B_1\left(\frac{1}{2}-\frac{1}{p+1},\frac{1}{2}\right)-B_{\frac{1}{2}}\left(\frac{1}{2}-\frac{1}{p+1},\frac{1}{2}\right)\right) \leq \tanh^{-1}\big(\frac{1}{\sqrt{2}}\big)\]
This has now given us the estimate for $p>1$, 
\begin{equation}
\int_{\Sigma_1} \IK^{p}(z,f) d\sigma_1    \geq   \big(4\pi(g-1)-4\big) + \frac{\pi^p}{4^p}  \frac{ 2\ell^{1-p}}{( \pi b)^p} 
\end{equation}
Here we have recalled from Corollary \ref{cor2} that the area of $\IA_1 \leq 4$.

\subsection{The approximation} Finally a word on the approximation of equation $u$ of (\ref{uxeqn}) by $v$ of (\ref{approxuxeqn}).  Notice that $0<v_x,u_x\leq 1$, that  $a\gg b$, and
\[ \int_{0}^{a} u_x = \int_{0}^{a} v_x =b \]
Indeed Corollary \ref{cor2} gives $a\geq \frac{1}{\ell}$, and $b$ is fixed and bounded.
Further  $ P[t] : t\mapsto \Big(1-\frac{1}{t^2}\Big)(t+\frac{1}{t})^{p-1}$, $t\in [0,1]$, 
is strictly increasing,  as is  $1/\lambda(x)$ on $[0,a/2]$. Hence $u_x$ is strictly decreasing on $[0,a/2]$ (and increasing on $[a/2,a]$). Fix $\beta<1$  a parameter.

Now $u_x\leq \beta$ on the interval $[b/\beta,a/2]$ (independent of $\ell$). We will choose $b/\beta\ll a$.  Then  
\[ P[t]=-\left(\frac{1}{t}\right)^{p+1} \left(1-(2-p) t^2 +\frac{1}{2} ((p-4) (p-1)) t^4+O\left(t^5\right)\right)\]
On the interval $[b/\beta,a/2]$ the term in brackets is uniformly bounded above and below, and these bounds both tend to $1$ as $\beta\to 0$.
Thus $u_x<\beta$ gives
\[ -(1+\epsilon)\left(\frac{1}{t}\right)^{p+1}\leq P[t]\leq -(1-\epsilon)\left(\frac{1}{t}\right)^{p+1} \]
and so the solutions $P[u_x]=c$, $u\approx c^{-\frac{1}{p+1}} = v_x$ and the distortion functions of the two solutions are $u_x+1/u_x$ and $v_x+1/v_x$ are now closely comparable. All this can be achieved for fixed $p$ when $\beta=\ell^s$, $s$ a small power (so $\epsilon \approx p\ell^{2s}$ and $\beta\approx \frac{\ell^{s/2}}{p}$).
The extra distortion we pick up in the interval $[0,b/\beta]$ is
\[ \int_{0}^{b/\beta} \big(\beta+\frac{1}{\beta}\big)^p  \approx \frac{b}{\beta^{p+1}} \]
 To retain our lower bounds no distortion,  we want this much less than $\ell^{1-p}$.  With $\beta=\ell^s$,  this is $\ell^{-s(p+1)} \ll \ell^{1-p}$,  so any $s\ll (p-1)(p+1)$ will eventually reveal our asymptotics as $\ell\to 0$.

\end{document}